\documentclass[10pt]{amsart}

\usepackage{amsmath,amscd,amssymb,amsfonts,amsthm,mathabx,hyperref}
\usepackage{xcolor} 

\usepackage{tcolorbox}
\usepackage[cmtip,matrix,arrow]{xy}
\usepackage{tikz}
\usepackage{tikz-cd}
\usepackage{graphicx} 
\usepackage{amsfonts,amsmath,amssymb,amsthm, mathtools}
\usepackage{cancel}
\usepackage{tikz-cd}
\usepackage{scalerel}
\usetikzlibrary{arrows.meta}
\usepackage{hyperref}
\usepackage{stackengine,wasysym}
\usepackage{tikz}
\usetikzlibrary{arrows.meta,decorations.markings}
\usepackage{wrapfig}
\usepackage{float}

\newcommand\reallywidetilde[1]{\ThisStyle{%
  \setbox0=\hbox{$\SavedStyle#1$}%
  \stackengine{-.1\LMpt}{$\SavedStyle#1$}{%
    \stretchto{\scaleto{\SavedStyle\m\kern.2mu\AC}{.5150\wd0}}{.6\ht0}%
  }{O}{c}{F}{T}{S}%
}}

\newtheorem{theorem}{Theorem}[section]
\newtheorem{proposition}[theorem]{Proposition}

\newtheorem{lemma}[theorem]{Lemma}
\newtheorem{rmk}[theorem]{Remark}

\theoremstyle{definition}
\newtheorem{definition}[theorem]{Definition}

\theoremstyle{remark}

\newtheorem{maintheorem}{Theorem}

\DeclareMathOperator{\Tr}{Tr}
\makeatletter
\def\paragraph{\@startsection{paragraph}{4}%
  \z@\z@{-\fontdimen2\font}%
  {\normalfont\bfseries}}
\makeatother

\newcommand{\ca}{\mathcal}

\newcommand{\C}{\mathbb{C}}

\renewcommand{\c}{\mathcal{C}}

\newcommand\lie[1]{\mathfrak{#1}}

\newcommand{\m}{\lie{m}}

\renewcommand{\u}{\lie{u}}
\newcommand{\on}{\operatorname}

\renewcommand{\ker}{ \on{ker}}

\newcommand{\ol}{\overline}

\newcommand{\p}{\partial}

\newcommand{\eeq}{\end{eqnarray*}}
\newcommand{\beq}{\begin{eqnarray*}}

\newcommand{\Diff}{\on{Diff}}

\newcommand{\PSU}{\on{PSU}}

\renewcommand{\S}{\ca{S}}

\newcommand{\ignore}[1]{}

\renewcommand{\subset}{\subseteq}

\newcommand{\DD}{\mathbb{D}}

\newcommand{\reals}{\mathbb{R}}

\begin{document}
\sloppy

\title[Teichmüller spaces with conical singularities and Virasoro coadjoint orbits]{Teichmüller spaces with conical singularities and Virasoro coadjoint orbits}
\author{Anton Alekseev}
\author{Rea Dalipi}
\author{Samson Shatashvili}

\begin{abstract}
We study Teichmüller spaces of hyperbolic \(0\)-metrics on a disk with ideal boundary and a marked point carrying a conical singularity with angle $\theta = 2\pi \beta$. For all $\beta>0$, we show that \(\operatorname{Teich}^{\beta}(\overline{\mathbb D})\cong S^1\backslash\operatorname{Diff}(S^1)\), and we construct  natural closed 2-forms $\omega_\beta$ induced by the Atiyah–Bott construction. For non-integer \(\beta\), the form $\omega_\beta$ is symplectic,  and one can identify \(\operatorname{Teich}^{\beta}(\overline{\mathbb D})\) with the corresponding elliptic Virasoro coadjoint orbit. 

For positive integers \(\beta=n\), $\omega_\beta$ is presymplectic, and it descends to the Kirillov-Kostant-Souriau symplectic form on the exceptional Virasoro orbit \(\mathrm{PSU}^{(n)}(1,1)\backslash\operatorname{Diff}(S^1)\). Infinitesimal isometries of a hyperbolic \(0\)-metric on a disk with a conical singularity with angle $\theta=2\pi n$ are given by the Lie algebra \(\mathfrak{sl}_2(\mathbb R)\) which integrates to the group $\mathrm{PSU}^{(n)}(1,1)$ on the boundary of the disk but not in the interior for $n>1$. It is remarkable that the presymplectic structure of the corresponding Teichmüller spaces recovers this group and descends to the corresponding Virasoro coadjoint orbit.

\end{abstract}
\maketitle 
\tableofcontents

\section{Introduction}
Teichmüller spaces play a central role in complex analysis and hyperbolic geometry, and have important applications in mathematical physics. Recent developments in theories of gravity have renewed interest in Teichmüller spaces of surfaces with boundary, which are closely related to Virasoro coadjoint orbits \cite{AS}.  In particular, Jackiw-Teitelboim (JT) gravity on the disk admits a boundary description governed by the Schwarzian action on the circle \cite{sss,stanford-witten, turiaci-mertens}. 
The classification of Virasoro coadjoint orbits is well established
\cite{latz-pank,witten}, with orbits falling into elliptic, hyperbolic,
parabolic and exceptional types.
The authors in \cite{AM} constructed natural symplectic structures on Teichmüller spaces of surfaces with ideal boundary and established their relationship with Virasoro geometry. In particular, they showed that the Teichmüller space of the disk is symplectomorphic to the first exceptional Virasoro coadjoint orbit.

In this paper, we extend this geometric picture by considering a disk with ideal boundary and a marked point carrying a conical singularity of total angle $2\pi\beta$, $\beta>0$. We describe the corresponding Teichmüller space and construct a natural closed two-form on it. For non-integer \(\beta>0\), the form \(\omega_\beta\) is symplectic, and the moment map identifies the Teichmüller space  with the elliptic Virasoro coadjoint orbit through the constant Hill potential \(T=\beta^2/4\). For \(\beta=n\in\mathbb Z_{>0}\), the moment map agrees with the natural projection onto the corresponding exceptional orbit, \(\omega_n\) is presymplectic and it is the pullback of its KKS symplectic form under this projection.

Related approaches to this correspondence have been developed using gauge theory. Valach and Youmans \cite{valach-youmans} showed that exceptional Virasoro coadjoint orbits arise through a Drinfeld--Sokolov reduction of $\mathrm{PSL}(2,\mathbb{R})$ BF theory on the punctured disk. Blau and Youmans \cite{blau-you} developed a gauge-theoretic description relating general Virasoro coadjoint orbits to hyperbolic geometries. They showed that infintesimally, isometries of the corresponding hyperbolic metrics are identified with the stabilizers of the coadjoint orbits.

In more detail, let $\overline{\mathbb{D}}$ be the closed disk and $p \in {\DD}$ a marked point in its interior. We consider hyperbolic metrics on $\ol\DD \setminus \{p\}$ with ideal boundary and 
 a simple singularity of order $\beta \in \mathbb{R}_{> 0}$  at $p$ (see Definition \ref{simple_singularity} below). Such singularities are called conical, with total angle $2\pi \beta$. Fixing $p$ to be at the origin, the model metric is  obtained by pulling back the Poincaré metric along the map $z \mapsto z^\beta$:
\begin{equation}\label{gbeta}
g_{\ol\DD}^\beta= \frac{4 \beta^2 |z|^{2(\beta-1)}}{(1-|z|^{2\beta})^2} \, |dz|^2.
\end{equation}
For positive integers $\beta=n$, the map $z\mapsto z^n$ is a
 branched covering.

Denote by $\rm Hyp^\beta $$(\ol \DD,p)$ the space of such metrics.  There is then a complex structure on 
$\ol \DD \setminus\{p\}$ induced by $g \in \rm Hyp^\beta $$(\ol \DD,p)$. It turns out that for simple singularities this complex structure smoothly extends to $\ol\DD$. We define the corresponding Teichmüller space as
\begin{equation*}
    \text{Teich}^\beta (\ol{\DD},p)=\text{Hyp}^\beta(\ol \DD,p)/\text{Diff}(\ol{\DD}, \partial \DD,p)
\end{equation*}
where $\rm Diff(\ol{\DD}, \partial \DD,p)$ is the group of orientation preserving diffeomorphisms fixing the boundary and the point $p$.

For $\beta=1$, the model metric is smooth at $p$. After
forgetting the marked point, one obtains the Teichmüller space
of the smooth disk with ideal boundary, studied in \cite{AM}.
The authors equipped this space with a natural
symplectic structure and  proved that it is isomorphic to the first exceptional Virasoro coadjoint orbit:
\begin{equation}
    {\rm Teich}(\ol{\DD})
    \cong
    \PSU(1,1) \backslash {\rm Diff}(S^1).
\end{equation}
Here ${\rm Diff}(S^1)$ denotes the group of orientation
preserving diffeomorphisms of the circle. Note that the
infinitesimal $\mathfrak{sl}_2(\mathbb{R})$ isometries of the
smooth Poincaré metric integrate to global $\PSU(1,1)$
Möbius isometries.

The goal of this paper is to extend this geometric picture to
elliptic and higher exceptional Virasoro coadjoint orbits. The exceptional orbits are of the form
\begin{equation*}
    \mathcal{O}^{(n)}
    =
    \PSU^{(n)}(1,1) \backslash {\rm Diff}(S^1),
\end{equation*}
where $\PSU^{(n)}(1,1)$ is the $n$-fold cover of the Möbius group
acting on the circle. One should note that for
$\beta=n \geq 2$, the infinitesimal isometries of the metric \eqref{gbeta}
again define an $\mathfrak{sl}_2(\mathbb{R})$-action on the
punctured disk. However, as it was shown in \cite{dal} the corresponding vector fields
are not all complete, and this action does not integrate to
a global $\PSU^{(n)}(1,1)$-action on the punctured disk.

Our first main result is the following.
\begin{maintheorem}
For every $\beta>0$, we have 
\begin{equation}\label{teich-virasoro}
    {\rm Teich}^{\beta}(\ol{\DD},p)
    \cong
    S^1 \backslash {\rm Diff}(S^1).
\end{equation}
This space carries a natural closed $2$-form $\omega_\beta$
induced by the Atiyah--Bott construction.
\end{maintheorem}

Note that the description \eqref{teich-virasoro} is independent of
$\beta$. However,  we show that $\beta$ reappears  in the
$2$-form $\omega_\beta$, and it allows us to
recover the distinction between elliptic and exceptional
Virasoro orbits. In particular, at positive integer values
of $\beta$, the 2-form $\omega_\beta$ becomes degenerate, and its null
directions determine the further quotient leading to an
exceptional Virasoro orbit. Our second main result makes this
relation precise.
\begin{maintheorem}
Let $\omega_\beta$ be the closed $2$-form on
${\rm Teich}^{\beta}(\ol{\DD},p)$ introduced above. 
\begin{enumerate}
    \item[(i)]
    If $\beta\notin\mathbb{Z}_{>0}$, the form $\omega_\beta$
    is symplectic and the moment map identifies
    ${\rm Teich}^{\beta}(\ol{\DD},p)$ with the corresponding
    elliptic Virasoro coadjoint orbit, equipped with its
    Kirillov--Kostant--Souriau form.

    \item[(ii)]
    If $\beta=n\in\mathbb{Z}_{>0}$, the form $\omega_n$
    is presymplectic,  and it satisfies
    \begin{equation*}
        \omega_n
        =
        \pi_n^*\omega_{\rm KKS}^{(n)},
    \end{equation*}
    where $\omega_{\rm KKS}^{(n)}$ is the
    Kirillov--Kostant--Souriau form on $\mathcal{O}^{(n)}$
    and
    \begin{equation*}
        \pi_n:
        {\rm Teich}^{n}(\ol{\DD},p)
        \cong S^1 \backslash {\rm Diff}(S^1)
        \longrightarrow
        \PSU^{(n)}(1,1) \backslash{\rm Diff}(S^1)
        =
        \mathcal{O}^{(n)}
    \end{equation*}
    is the natural projection. 
\end{enumerate}
\end{maintheorem}
Remarkably, despite the nonintegrability of the infinitesimal
isometries for $\beta=n>1$ discussed above, this result
recovers both the group $\PSU^{(n)}(1,1)$ and the dependence
on $n$, neither of which is visible in the homogeneous-space
description ${\rm Teich}^{n}(\ol{\DD},p) \cong S^1 \backslash {\rm Diff}(S^1)$.

The structure of the paper is as follows. In Section 2, we recall the relevant properties of Virasoro coadjoint orbits, with particular attention to the exceptional orbits. In Section 3, we define Teichmüller spaces of hyperbolic  metrics with singularities and establish their homogeneous-space descriptions. Section 4 contains the construction of the two-forms $\omega_\beta$ and the proof of the main result. Appendix A recalls the uniformization results used in Section 3.

\vspace{4mm}

\paragraph{Acknowledgements}
We would like to thank Matthias Blau, Eckhard Meinrenken, Julian Sonner, Marc Troyanov, Herman Verlinde, Yilin Wang and Donald Youmans  for sharing their ideas and for inspiring discussions on the topic. R.D. greatly appreciates the hospitality of Hamilton Mathematics Institute at Trinity College Dublin where part of this work was completed. 

Research of A.A. and R.D. was supported in part by the grant 236683
of the Swiss National Science Foundation (SNSF). Research of A.A. was supported in part by the grant 10009199 of the SNSF. Research of A.A. and S.S. was supported in part by the award of the Simons Foundation to the Hamilton Mathematics Institute of the Trinity College Dublin under the program “Targeted Grants to Institutes”, and research of S.S. was supported by Hamilton Mathematics Institute at Trinity College Dublin.

\section{Virasoro coadjoint orbits}

\subsection{Virasoro coadjoint action and coadjoint orbits}

In this Section, we collect some standard facts about the central extension of the group of diffeomorphisms of the circle, of its coadjoint Virasoro action, and of the corresponding coadjoint orbits \cite{bal-feh-pal, latz-pank, Kirillov}.

\subsubsection{The Bott-Virasoro central extension}

The group ${\rm Diff}(S^1)$ of orientation preserving diffeomorphisms can be represented as
$$
{\rm Diff}(S^1) = \{ \phi: S^1 \to S^1; \phi'(x)>0, \int_{S^1} \phi'(x) dx = 2\pi\}.
$$
Its fundamental group is given by $\pi_1({\rm Diff}(S^1)) \cong \mathbb{Z}$, and its universal cover can be identified with the group of quasi-periodic diffeomorphisms of $\mathbb{R}$ with period $2\pi$:
$$
\widetilde{\rm Diff}(S^1)=\{ f: \mathbb{R} \to \mathbb{R}; f'(x) >0, f(x+2\pi) = f(x) + 2\pi\}.
$$
A rigid rotation by angle $\theta$ is represented by an element 
$$
r_\theta(x) = x + \theta.
$$
Rigid rotations $r_{2\pi k}$ by multiples of $2\pi$ are central elements of $\widetilde{\rm Diff}(S^1)$, and we have
$$
{\rm Diff}(S^1) \cong \widetilde{\rm Diff}(S^1)/2\pi \mathbb{Z} \cong 2\pi \mathbb{Z} \backslash \widetilde{\rm Diff}(S^1).
$$
If one uses the map $x \mapsto z=e^{ix}$ to identify $S^1$ with the unit circle on the complex plane, the action of ${\rm Diff}(S^1)$ acquires the form
$$
\phi: z=e^{ix} \mapsto w=e^{i\phi(x)}.
$$
The group ${\rm Diff}(S^1)$ has a subgroup isomorphic to ${\rm PSU}(1,1)$ acting by Möbius transformations
\begin{equation}      \label{eq:Möbius}
z \mapsto \frac{az +b}{\bar{b} z + \bar{a}},
\end{equation}
where $|a|^2-|b|^2=1$. It also has subgroups isomorphic to $n$-fold covers ${\rm PSU}^{(n)}(1,1)$ of ${\rm PSU}(1,1)$. Their actions are defined by formulas
$$
z \mapsto \left(\frac{az^n +b}{\bar{b} z^n + \bar{a}}\right)^{\frac{1}{n}}.
$$
Each branch of the $n$-th root function defines an element of ${\rm PSU}^{(n)}(1,1)$ covering the elements of ${\rm PSU}(1,1)$ defined by $a$ and $b$.

Both ${\rm Diff}(S^1)$ and $\widetilde{\rm Diff}(S^1)$ admit  central extensions defined by the Bott-Virasoro 2-cocycle
\begin{equation}
c(f_1, f_2) = \frac{1}{2} \int_{S^1} \log(f'_1(f_2(x))) \, \frac{d}{dx} \, \log(f'_2(x))  \, dx.
\end{equation}
Here one can replace $f: \mathbb{R} \to \mathbb{R}$ by $\phi: S^1 \to S^1$.
Note that for every $f$ one has $c(r_\theta, f)=c(f, r_\theta)=0$ for all rigid rotations of the circle.

\subsubsection{The Virasoro coadjoint action}

The Virasoro coadjoint action is defined on the space of second order Hill operators on the circle
$$
{\rm Hill}(S^1) = \left\{ \partial_x^2 + T(x); T(x+2\pi) = T(x)\right\},
$$
and it can be represented as follows
\begin{equation}         \label{eq:Virasoro action}
\begin{array}{llll}
f^{-1}:  &\partial_x^2 + T(x) & \mapsto  &f'(x)^{3/2}\left(\frac{\partial^2}{\partial f(x)^2} + T(f(x))\right)f'(x)^{1/2} \\
 & & = & \partial_x^2 + \left(T(f(x))f'(x)^2 + \frac{1}{2} \mathcal{S}(f)\right),
 \end{array}
\end{equation}
where $\mathcal{S}(f)$ is the Schwarzian derivative
\begin{equation}\label{schwarzian}
    \mathcal{S}(f) = \frac{f'''}{f'} - \frac{3}{2} \frac{f''^2}{f'^2}.
\end{equation}
Note that the right hand side of \eqref{eq:Virasoro action} defines a right action of the element $f$, and a left action of the element $f^{-1}$.

For a Hill potential $T(x)$,  the  Virasoro coadjoint orbit $\mathcal{O}_T$ comes equipped with the canonical projection
$$
\nu_T: {\rm Diff}(S^1)\to \mathcal{O}_T = {\rm Stab}_T\backslash  {\rm Diff}^+(S^1),
$$
where ${\rm Stab}_T$ is the stabilizer of the Hill operator $\partial_x^2+T$ under the Virasoro coadjoint action. 

We will focus our attention on Virasoro coadjoint orbits passing through constant Hill potentials $T(x)={\rm const}$ (for the full classification see \cite{latz-pank}).
For $T<0$, the corresponding orbits are called hyperbolic, for $T>0$ but $T\notin \{ \frac{n^2}{4}; n \in \mathbb{N}\}$ they are called elliptic, and in both cases
we have ${\rm Stab}_T = S^1$ realized by rigid rotations of the circle.

Virasoro orbits with $T\in \{ \frac{n^2}{4}; n \in \mathbb{N}\}$ are called exceptional. We will use the notation $\mathcal{O}^{(n)}$ for the orbit passing through $T=\frac{n^2}{4}$. Its stabilizer subgroup is given by
$$
{\rm Stab}^{(n)} = {\rm PSU}^{(n)}(1,1) \subset {\rm Diff}(S^1).
$$
In particular, the orbit $\mathcal{O}^{(1)}$ is called the Teichmüller orbit, and ${\rm Stab}^{(1)} = {\rm PSU}(1,1)$ is the Möbius group.

\subsubsection{KKS symplectic forms}

A coadjoint orbit $\mathcal{O}_T$ carries a canonical Kirillov-Kostant-Souriau (KKS) symplectic form $\omega_T$ such that
\begin{equation} \label{eq: KKS}
\nu_T^* \, \omega_T = \frac{1}{4} \int \delta(\log(f'(x)) \wedge \delta(\log(f'(x))') dx + 
\int T(f(x)) \, \delta f'(x) \wedge \delta f(x) dx.
\end{equation}
Under this symplectic structure, the Hamiltonian of rigid rotations is of the form
\begin{equation}        \label{eq: Schwarzian}
H= \int \left( T(f(x))f'(x)^2 + \frac{1}{2} \mathcal{S}(f) \right) dx.
\end{equation}

\subsection{Exceptional Virasoro orbits}

In this Section, we focus our attention on the structure of exceptional Virasoro coadjoint orbits.

\subsubsection{${\rm PSU}^{(n)}(1,1)$ action on the circle}

Recall that the action of ${\rm PSU}^{(n)}(1,1)$ is defined using the branches of the $n$th  root function.
In more detail,  let 
$$
\mu_n: z \mapsto z^n, \hskip 0.3cm \Pi_n: {\rm PSU}^{(n)}(1,1) \to {\rm PSU}(1,1)
$$
be the $n$th power function and the natural projection of ${\rm PSU}^{(n)}(1,1)$ to ${\rm PSU}(1,1)$. The following simple proposition relates their actions on the circle.
\begin{proposition}     \label{prop:mu_n pi_n}
For all $z \in S^1, g \in {\rm PSU}^{(n)}(1,1)$ we have
$$
\mu_n(g \cdot z) = \Pi_n(g) \cdot \mu_n(z).
$$
\end{proposition}

\begin{proof}
The proof is by a direct computation:
$$
\mu_n(g \cdot z) = \mu_n\left(\frac{az^n +b}{\bar{b} z^n + \bar{a}}\right)^{\frac{1}{n}} =
\frac{az^n +b}{\bar{b} z^n + \bar{a}} = \Pi_n(g) \cdot z^n= \Pi_n(g) \cdot \mu_n(z).
$$
\end{proof}

\subsubsection{The spaces ${\rm Diff}^{(n)}(S^1)$ and $\widetilde{\rm Diff}^{(n)}(S^1)$.}

For future use, we introduce the spaces
$$
\begin{array}{lll}
\widetilde{\rm Diff}^{(n)}(S^1) & = & \{ F: \mathbb{R} \to \mathbb{R}; F'(x) >0, F(x+2\pi) = F(x) + 2\pi n \}, \\
{\rm Diff}^{(n)}(S^1) & = & \{ \phi: S^1 \to S^1; \phi'(x) >0, \int_{S^1} \phi'(x) dx = 2\pi n \},
\end{array}
$$
and the maps
$$
m_n: \mathbb{R} \to \mathbb{R}, \hskip 0.3cm m_n(x) = nx.
$$
Note that $m_n \in \widetilde{\rm Diff}^{(n)}(S^1)$ and $\mu_n \in {\rm Diff}^{(n)}(S^1)$.
Properties of the spaces $\widetilde{\rm Diff}^{(n)}(S^1)$ and ${\rm Diff}^{(n)}(S^1)$ are collected in the following proposition:
\begin{proposition}\label{circle-n-folds}
The map $q_n: f \mapsto m_n \circ f$ establishes an isomorphism $\widetilde{\rm Diff}^{(n)}(S^1) \cong \widetilde{\rm Diff}(S^1)$.
Furthermore, for all $f_1, f_2 \in \widetilde{\rm Diff}(S^1)$ we have
$$
q_n(f_1 \circ f_2) = q_n(f_1) \circ f_2.
$$
The map $q_n$ descends to the $n$-fold cover map
$$
p_n: {\rm Diff}(S^1) \to {\rm Diff}^{(n)}(S^1) \cong \mathbb{Z}_n\backslash {\rm Diff}(S^1)
$$
given by $p_n(\phi) = \mu_n \circ \phi$. For all $\phi, \phi_1, \phi_2 \in  {\rm Diff}(S^1), g \in  {\rm PSU}^{(n)}(1,1)$ we have
$$
p_n(\phi_1 \circ \phi_2) = p_n(\phi_1) \circ \phi_2, \hskip 0.3cm p_n( g \circ \phi) = \Pi_n(g) \circ p_n(\phi).
$$
\end{proposition}

\begin{proof}
First, note that the map $q_n$ does map $\widetilde{\rm Diff}(S^1)$ to $\widetilde{\rm Diff}^{(n)}(S^1)$. Indeed,
let $F(x) = q_n(f)(x)=m_n((f)(x)) = nf(x)$. Then,
$$
F(x+2\pi) = nf(x+2\pi) = n(f(x) + 2\pi) = nf(x) + 2\pi n = F(x) + 2\pi n,
$$
as required. The map $q_n$ is invertible, its inverse is given by
$$
q_n^{-1}(F)(x) = \frac{1}{n} F(x).
$$
Hence, $q_n$ is an isomorphism between $\widetilde{\rm Diff}(S^1)$ and $\widetilde{\rm Diff}^{(n)}(S^1)$. 
Next, we compute
$$
q_n(f_1 \circ f_2)(x) = nf_1(f_2(x))=q_n(f_1)(f_2(x)) = (q_n(f_1) \circ f_2)(x).
$$

A similar computation shows that $p_n(\phi_1 \circ \phi_2) = p_n(\phi_1) \circ \phi_2$. For the action of ${\rm PSU}^{(n)}(1,1)$, we have
$$
p_n( g \circ \phi) = \mu_n \circ g \circ \phi = \Pi_n(g) \circ \mu_n \circ \phi = \Pi_n(g) \circ p_n(\phi),
$$
where we have used Proposition \ref{prop:mu_n pi_n}. Finally, the preimage of  $\mu_n \in {\rm Diff}^{(n)}(S^1)$ under the map $p_n$ is the subgroup 
$\mathbb{Z}_n \subset {\rm Diff}(S^1)$ given by
rigid rotations $r_{2\pi k/n}$ for $k=0, \dots, n-1$ showing that ${\rm Diff}^{(n)}(S^1) \cong \mathbb{Z}_n\backslash {\rm Diff}(S^1)$.
\end{proof}

\subsubsection{Structure of $\mathcal{O}^{(n)}$}

The following proposition explains the relation between exceptional orbits and spaces ${\rm Diff}^{(n)}(S^1)$:
\begin{proposition}
For the exceptional coadjoint orbit $\mathcal{O}^{(n)}$ we have
$$
\mathcal{O}^{(n)} \cong  {\rm PSU}(1,1) \backslash {\rm Diff}^{(n)}(S^1).
$$
Under the projection $\sigma_n: {\rm Diff}^{(n)}(S^1) \to \mathcal{O}^{(n)}$ the symplectic form reads
$$
\sigma_n^* \omega^{(n)}_{\rm KKS} = \frac{1}{4} \int_{S^1} \left(\delta(\log(F'(x)) \wedge \delta(\log(F'(x))')  + \delta F'(x) \wedge \delta F(x)\right) \, dx,
$$
and the moment map of rigid rotations acquires the form
$$
H = \int_{S^1} \left( \frac{1}{4} F'^2(x) + \frac{1}{2} \mathcal{S}(F)\right) dx.
$$
\end{proposition}

\begin{proof}
Recall that ${\rm Diff}^{(n)}(S^1) \cong \mathbb{Z}_n \backslash {\rm Diff}(S^1)$, where $\mathbb{Z}_n$ acts by rigid rotations by angles $2\pi k/n$ which is exactly the action of the center of  the subgroup ${\rm PSU}^{(n)}(1,1) \subset {\rm Diff}(S^1)$. Hence,
$$
O^{(n)} \cong {\rm PSU}^{(n)}(1,1) \backslash {\rm Diff}(S^1) \cong ({\rm PSU}^{(n)}(1,1)/\mathbb{Z}_n) \backslash {\rm Diff}^{(n)}(S^1) =
{\rm PSU}(1,1) \backslash {\rm Diff}^{(n)}(S^1),
$$
as required. The formulas for the symplectic form and for the Hamiltonian of rigid rotations are obtained by substitution of $F(x) = nf(x)$ in equations \eqref{eq: KKS} and \eqref{eq: Schwarzian}.
\end{proof}

\section{Teichmüller spaces of a disk}
\label{Gravity_theory}

\subsection{Hyperbolic structures and Teichmüller spaces}
In this Section, we recall the notions of hyperbolic $0$-metrics and of (infinite dimensional) Teichmüller spaces associated to surfaces with boundary.

\subsubsection{$0$-metrics}
The model for hyperbolic structures on surfaces with boundary is the closed Poincar\'{e} disk 
\[ \ol{\DD}=\{z\in \C|\ |z|\leq1\},\] 
equipped with the standard action of $G=\PSU(1,1)$ by M\"obius transformations \eqref{eq:Möbius}. Recall that it carries the standard Poincaré metric
\begin{equation}      \label{eq:Poincaré metric}
g_\mathbb{D} = \frac{4|dz|^2}{(1-|z|^2)^2},
\end{equation}
and that $G$ is the group of isometries of $g_\mathbb{D}$.

\begin{definition}\label{def:hyperbolic}
	A \emph{hyperbolic structure} on an oriented surface $\Sigma$ with boundary  $\p\Sigma$ is an oriented atlas with  $\ol{\DD}$-valued charts, with constant transition maps given by elements of $G$. The space of all hyperbolic structures on $\Sigma$  will be denoted 
	$\on{Hyp}(\Sigma)$. 
\end{definition}

For $U \subset \Sigma$ a chart on $\Sigma$, a local map (a developing map) $\sigma: U \to \ol{\DD}$ induces a Riemannian metric on $U \cap \mathrm{int}(\Sigma)$ by pulling back the Poincar\'e metric $g_\mathbb{D}$. Riemannian metrics defined on different charts coincide on overlaps and give rise to a hyperbolic Riemannian metric on $\mathrm{int}(\Sigma)$. Boundary components become ideal boundaries, that is, boundaries at infinity arising from the compactification of the hyperbolic metric. Such metrics exhibit a specific boundary behavior, captured by the concept of a \textit{$0$-metric}.

\begin{definition}\label{0-metric}
Let $\Sigma$ be a smooth manifold with boundary. A Riemannian metric $g$ on the interior $\mathrm{int}(\Sigma)$ is called a \emph{0-metric} if for every smooth boundary defining function $\varrho$ on $\Sigma$ the tensor $\varrho^2 g$ smoothly extends  to a  metric on $\Sigma$.
A 0-metric is \emph{hyperbolic} if its Gauss curvature satisfies $K_{\mathsf{g}} = -1$ on $\mathrm{int}(\Sigma)$.\footnote{We recall that a boundary defining function is a smooth nonnegative funciton $\rho$ such that $\rho^{-1}(0)=\partial \Sigma$ and $d \rho \neq 0$ along the boundary $\partial \Sigma$.}
\end{definition}

Clearly, the Poincaré metric   \eqref{eq:Poincaré metric}
is a 0-metric on $\mathbb{D}$. Indeed, 
for a boundary defining function $\varrho=1-|z|^2$ we have
$$
\varrho^2 g_\mathbb{D} = {4 |dz|^2}
$$
which extends  smoothly to the boundary of $\ol{\DD}$.

\subsubsection{Teichmüller spaces}
For a surface with boundary $\Sigma, \partial \Sigma \neq \emptyset$, we denote by ${\rm Hyp}(\Sigma)$ the space of hyperbolic $0$-metrics on $\Sigma$. 
 This space carries a natural action of the group ${\rm Diff}(\Sigma)$ of orientation preserving diffeomorphisms of $\Sigma$, and of its connected component ${\rm Diff}_0(\Sigma)$.
 
 We define the corresponding Teichmüller space
 \begin{equation}
     {\rm Teich}(\Sigma) ={\rm Hyp}(\Sigma)/{\rm Diff}_0(\Sigma; \partial \Sigma),
 \end{equation}
 where 
 $$
{\rm Diff}_0(\Sigma; \partial \Sigma) = \{ \varphi \in {\rm Diff}(\Sigma);
\varphi|_{\partial \Sigma} = {\rm Id}\}_0
 $$
 is the connected component of the group ${\rm Diff}(\Sigma; \partial \Sigma)$ of orientation preserving diffeomorphisms of $\Sigma$ which restrict to the identity on the boundary.
Note that in the case of $\Sigma = \ol{\DD}$, the groups ${\rm Diff}(\ol{\DD})$ and ${\rm Diff}(\ol{\DD}; \partial \ol{\DD})$ are connected, and we omit the subscripts.
We recall the following elementary fact:
\begin{proposition}
    Assume that $\partial \Sigma \cong S^1$. Then, for surfaces of positive genus
    $$
{\rm Diff}_0(\Sigma)/{\rm Diff}_0(\Sigma; \partial \Sigma) \cong \widetilde{\rm Diff}(S^1),
    $$
    and for the disk
    \begin{equation} \label{eq:disk quotient}
{\rm Diff}(\ol{\DD})/{\rm Diff}(\ol{\DD}; \partial \ol{\DD}) \cong {\rm Diff}(S^1).
    \end{equation}
\end{proposition}

\begin{proof}
Denote $\mathcal{G}={\rm Diff}_0(\Sigma), \mathcal{H}={\rm Diff}_0(\Sigma; \partial \Sigma)$ and
$$
\mathcal{F} = {\rm Diff}_0(\Sigma)/{\rm Diff}_0(\Sigma; \partial \Sigma).
$$
Restriction to the boundary gives a surjective homomorphism
\[
r:\mathcal{G}\longrightarrow \operatorname{Diff}(S^1),
\qquad r(\phi)=\phi|_{\partial\Sigma},
\]
which admits local sections by isotopy extension. Since
\(\mathcal{H}=(\ker r)_0\), this induces a connected covering
\[
\overline r:\mathcal{F}=\mathcal{G}/\mathcal{H}\longrightarrow \operatorname{Diff}(S^1).
\]
Recall the exact sequence in homotopy groups
$$
\dots \to \pi_1(\mathcal{H}) \to \pi_1(\mathcal{G}) \to \pi_1(\mathcal{F}) \to \pi_0(\mathcal{H}) \to \dots
$$
Since $\mathcal{H}={\rm Diff}_0(\Sigma; \partial \Sigma) $ is contractible we have $\pi_0(\mathcal{H}) = \pi_1(\mathcal{H}) = 1$ and $\pi_1(\mathcal{F}) \cong \pi_1(\mathcal{G})$.
For $\Sigma$ of positive genus, $\chi(\Sigma)=1-2g <0$, so
$\mathcal{G}$ is contractible \cite{homotopy}.
Hence, $\pi_1(\mathcal{F})=1$, and we conclude that $\ol r$ is the universal cover,
$$
\mathcal{F} \cong \widetilde{\rm Diff}(S^1).
$$
 For $\Sigma = \ol{\DD}$, we have $\pi_1(\mathcal{G}) \cong \mathbb{Z}$ realized by rigid rotations by multiples of $2\pi$. Hence,
$\pi_1(\mathcal{F}) \cong \mathbb{Z}$ 
and
$$
\mathcal{F} \cong {\rm Diff}(S^1).
$$
\end{proof}

\subsubsection{Teichmüller space of the disk}
Recall from \cite{AM} the following simple calculation of the Teichmüller space associated to the disk $\Sigma = \ol{\DD}$.
\begin{proposition}
    We have
    $$
{\rm Teich}(\ol{\DD}) \cong {\rm PSU}(1,1)\backslash {\rm Diff}(S^1).
    $$
\end{proposition}

Note that ${\rm Teich}(\ol{\DD})$ is diffeomorphic to  the Virasoro Teichmüller orbit $\mathcal{O}^{(1)}$.

\begin{proof}
Since the group of isometries of $g_\mathbb{D}$ is the Möbius group
${\rm PSU}(1,1)$, we have
    $$
{\rm Hyp}(\ol \DD) = g_\mathbb{D} \cdot {\rm Diff}(\ol{\DD}) \cong
{\rm PSU}(1,1)\backslash {\rm Diff}(\ol{\DD}).
    $$
Then, using \eqref{eq:disk quotient} we obtain
    $$
{\rm Teich}(\ol{\DD}) = {\rm Hyp}(\ol{\DD})/{\rm Diff}(\ol{\DD}, \partial \ol{\DD}) \cong {\rm PSU}(1,1)\backslash {\rm Diff}(\ol{\DD})/{\rm Diff}(\ol{\DD}, \partial \ol{\DD}) \cong
{\rm PSU}(1,1)\backslash {\rm Diff}(S^1),
    $$ 
    as required.
\end{proof}
Let $J$ denote the standard complex structure on $\ol \DD$ defined by the holomorphic coordinate $z$.
We recall that the stabilizer group of $J$ is also the Möbius group ${\rm PSU}(1,1)$, and that the space of complex structures on the disk $\mathcal{J}(\ol \DD)$ is given by
$$
\mathcal{J}(\ol\DD) = J \cdot {\rm Diff}(\ol{\DD}) \cong
{\rm PSU}(1,1)\backslash {\rm Diff}(\ol{\DD}) \cong {\rm Hyp}(\ol{\DD}).
$$

We will also need Teichmüller spaces corresponding to surfaces with marked points. For the disk, we introduce
the space of metrics 
$$
{\rm Hyp}^{(1)}(\ol \DD) = {\rm Hyp}(\ol \DD) \times \ol \DD,
$$
where an element $(g, p) \in {\rm Hyp}^{(1)}(\ol \DD)$ is a pair consisting of a hyperbolic $0$-metric $g$ on $\ol \DD$ and a point $p \in \mathbb{D}$. The Teichmüller space is defined as the quotient
$$
{\rm Teich}^{(1)}(\ol{\DD}) =  {\rm Hyp}^{(1)}(\ol \DD)/{\rm Diff}(\ol{\DD}; \partial \ol{\DD})
$$
under the diagonal action of the diffeomorphism group on metrics and on the disk.

\begin{proposition}       \label{prop:Teich^1}
    We have
    $$
{\rm Teich}^{(1)}(\ol{\DD}) \cong S^1 \backslash {\rm Diff}(S^1).
    $$
\end{proposition}

\begin{proof}
    As before, we use the fact that ${\rm Diff}(\ol{\DD})$ acts transitively on ${\rm Hyp}(\ol \DD)$. Furthermore, the stabilizer of the standard metric $g_\mathbb{D}$ is the Möbius group which acts transitively on the disk with stabilizers isomorphic to $S^1$. Hence,
    $$
{\rm Hyp}^{(1)}(\ol \DD) \cong S^1 \backslash {\rm Diff}(\ol{\DD}).
    $$
    Then,
    $$
{\rm Teich}^{(1)}(\ol{\DD}) = S^1 \backslash {\rm Diff}(\ol{\DD})/{\rm Diff}(\ol{\DD}, \partial \ol{\DD}) \cong
S^1 \backslash {\rm Diff}(S^1),
    $$
    as required.
\end{proof}

Note that we can again consider complex structures $\mathcal{J}^{(1)}(\ol\DD)$ on $\ol{\DD}$ with a marked point inside the disk. Then,
$$
\mathcal{J}^{(1)}(\ol\DD) = J \cdot {\Diff}(\ol{\DD}) \cong
S^1 \backslash {\rm Diff}(\ol{\DD}) \cong {\rm Hyp}^{(1)}(\ol{\DD}).
$$

\subsection{Teichmüller spaces with conical singularities}

In this Section, we discuss Teichmüller spaces associated to a disk with one conical singularity.

\subsubsection{Conical singularities}
Let $\Sigma$ be an oriented 2-manifold equipped with a complex structure $I$.

\begin{definition}
    A Riemannian metric $g$ on $\Sigma$ is conformal with respect to the complex structure $I$ if it can be presented in the form
    $$
g = e^{2v} |dz|^2,
    $$
    where $v$ is a real valued function, and $z$ is a local holomorphic coordinate on $\Sigma$.
\end{definition}
We will be interested in conformal metrics with singularities of a special type:

\begin{definition}\label{simple_singularity}
   A conformal metric $g$ on $\Sigma$ admits a \emph{simple singularity of order $\beta \in \reals $} at $p \in \Sigma$ if it can be locally written as
\[
g = e^{2u} |z|^{2(\beta-1)} |dz|^2,
\]
where $z$ is a local holomorphic coordinate with $z(p)=0$,
and  $u \in L_\text{loc}^1, \Delta u \in L_\text{loc}^1$. The condition $\Delta u \in L_{\text{loc}}^1$ is understood in the sense of distributions on a full neighborhood of $p$, including $p$ itself. If $\beta > 0$,  a simple singularity  is said to be \emph{conical} with cone angle $\theta=2\pi \beta$. \footnote{Our convention differs from the classical convention:
a singularity of order $\beta$ here has order $\beta-1$
in \cite{hulin-troyanov}.}
\end{definition}

Note that $\beta = 1$ corresponds to the angle $\theta = 2\pi$, and the factor $z^{2(\beta-1)}$ disappears. In the hyperbolic case considered below, the singularity is then removable (see Lemma \ref{lem:cone-local-regularity}) and the metric extends smoothly across $p$. For $0 < \beta < 1$, the angle $\theta = 2\pi \beta < 2\pi$, and the factor $z^{2(\beta -1)}$ diverges at $z=0$. For $\beta >1$, the angle $\theta = 2\pi \beta > 2 \pi$, and the factor $z^{2(\beta -1)}$ vanishes at $z=0$.

We will now focus on the case of $\Sigma = \ol{\DD}$ with one conical singularity. The following theorem summarizes several results of Troyanov et al. \cite{troyanov, hulin-troyanov}. Its proof, together with other relevant results from the theory of metrics with simple singularities, is given in Appendix \ref{app-conemetrics}.

\begin{theorem}       \label{thm:Troyanov}
    For a triple of a complex structure $I$ on $\ol\DD$, a point $p \in \DD$, and $\beta >0$, there is a unique hyperbolic conformal $0$-metric $g$ on $\ol\DD$ with a conical singularity of order $\beta$ at $p$.
\end{theorem}

As an example, consider the standard complex structure $J$ on the disk, and the point $p$ at the origin. Then, the locally defined holomorphic map $\phi_\beta: z \mapsto z^\beta$ gives rise to a conformal hyperbolic metric
$$
g_\DD^\beta = \phi_\beta^* g_\DD = \frac{4 \beta^2 |z|^{2(\beta-1)}}{(1-|z|^{2\beta})^2} \, |dz|^2.
$$
This metric has a conical singularity of order $\beta$, and the corresponding function
$$
u(z) = \log\left( \frac{2 \beta}{1 -|z|^{2\beta}}\right)
$$
is continous at $z=0$. Furthermore, observe that $g_\DD^\beta$ is a $0$-metric. Indeed, for the boundary defining function $\varrho=1-|z|$, we have
$$
\varrho^2  g_\DD^\beta =
\frac{4 \beta^2 |z|^{2(\beta-1)}(1-|z|)^2}{(1-|z|^{2\beta})^2} \, |dz|^2
$$
which is a regular Riemannian metric near the boundary $|z|=1$.

\subsubsection{Teichmüller spaces on the disk with conical singularities}

For future use, we introduce the following notation: let ${\rm Hyp}(\ol{\DD} \backslash p)$ be the space of hyperbolic $0$-metrics on the complement of the point $p \in \ol{\DD}$ which are conformal under some complex structure $I$ on $\ol{\DD}$ and have a conical singularity at $p$.
We denote the order of that singularity by ${\rm ord}(p)$.

We extend the definition of ${\rm Hyp}^{(1)}(\mathbb{D})$ to the space of hyperbolic metrics on the disk with a marked point having a conical singularity of order $\beta$:
$$
{\rm Hyp}^{\beta}(\ol{\DD})=
\{ (g, p); p \in \ol{\DD}, g \in {\rm Hyp}(\ol{\DD} \backslash p),
{\rm ord}(p) = \beta\}.
$$
Note that for $\beta=1$ the metric has a removable singularity (see Lemma \ref{lem:cone-local-regularity}), but the marked point data is still retained in our definition, and we recover $\rm Hyp^{(1)}(\mathbb{D})$.

Furthermore, we can consider hyperbolic metrics with a conical singularity at a given point $p$:
$$
{\rm Hyp}^{\beta}(\ol{\DD}, p)=
\{ g \in {\rm Hyp}(\ol{\DD} \backslash p),
{\rm ord}(p) = \beta\}.
$$

We also define the subgroup of the diffeomorphism group 
${\rm Diff}(\ol{\DD}; p) \subset {\rm Diff}(\ol{\DD})$ preserving $p$.

\begin{proposition}\label{hypbeta}
    We have
    $$
{\rm Hyp}^{\beta}(\ol\DD) \cong S^1 \backslash {\rm Diff}(\ol{\DD}),
\hskip 0.3cm 
{\rm Hyp}^{\beta}(\ol\DD, p) \cong S^1 \backslash {\rm Diff}(\ol{\DD}; p).
    $$
\end{proposition}

\begin{proof}
    By Theorem \ref{thm:Troyanov}, metrics with conical singularities on $\ol\DD$ are in bijection with complex structures with a marked point. Hence, we have an isomorphism
    $$
{\rm Hyp}^{\beta}(\ol\DD) \cong \mathcal{J}^{(1)}(\ol\DD) \cong 
S^1 \backslash {\rm Diff}(\ol{\DD}),
    $$
    as required. A similar argument applies to ${\rm Hyp}^{\beta}(\ol\DD, p)$.
\end{proof}

We define Teichmüller spaces with conical singularities as
$$
{\rm Teich}^\beta(\ol{\DD}) = \frac{{\rm Hyp}^\beta(\ol\DD)}{{\rm Diff}(\ol{\DD}; \partial\ol{\DD})}, \hskip 0.3cm
{\rm Teich}^\beta(\ol{\DD}, p) = \frac{{\rm Hyp}^\beta(\ol\DD, p)}{{\rm Diff}(\ol{\DD}; \partial\ol{\DD}, p)},
$$
where ${\rm Diff}(\ol{\DD}; \partial\ol{\DD}, p)$ is the subgroup of 
${\rm Diff}(\ol{\DD})$ preserving both the boundary pointwise and $p$. 
Note that since both 
${\rm Diff}(\ol{\DD}; \partial\ol{\DD})$ and ${\rm Diff}(\ol{\DD}; \partial\ol{\DD}, p)$ are connected, we omit the subscript $0$ indicating  the connected component of the identity.

\begin{proposition}
    For all $p \in \DD$ and $\beta >0$, we have
    $$
{\rm Teich}^\beta(\ol{\DD}, p) \cong {\rm Teich}^\beta(\ol{\DD}) \cong S^1 \backslash {\rm Diff}(S^1).
    $$
\end{proposition}

\begin{proof}
    For ${\rm Teich}^\beta(\ol{\DD}) \cong S^1 \backslash {\rm Diff}(S^1)$,  the proof repeats verbatim the one of Proposition \ref{prop:Teich^1}. For ${\rm Teich}^\beta(\ol{\DD}, p) \cong {\rm Teich}^\beta(\ol{\DD})$, note that ${\rm Diff}(\ol{\DD}; \partial\ol{\DD})$ acts transitively on the disk, so the conical singularity can always be moved to $p$, and the subgroup of ${\rm Diff}(\ol{\DD}; \partial\ol{\DD})$ preserving $p$ is exactly ${\rm Diff}(\ol{\DD}; \partial\ol{\DD}, p)$.
\end{proof}

\subsubsection{The case of $\beta=n \in \mathbb{Z}_{\geq 1}$}

In this paragraph, we consider the situation of integer values of $\beta$ in more detail. Fix $p$ at the origin of $\DD$, and let $u_n: \ol{\DD} \to \ol{\DD}$ be the map $u_n: z \mapsto z^n$,
and denote

\[
\operatorname{Diff}^{(n)}(\overline{\mathbb D},p)
:=
\{u_n\circ\phi:
\phi\in\operatorname{Diff}(\overline{\mathbb D},p)\}.
\]
%


\begin{proposition}
 Similarly as in Proposition \ref{circle-n-folds}, there is a map
\[
v_n:\mathrm{Diff}(\ol \DD,p)
\longrightarrow \mathrm{Diff}^{(n)}(\ol \DD,p),
\qquad
v_n(\varphi)=u_n\circ\varphi,
\]
inducing a bijection
\[
\mathbb Z_n\backslash\mathrm{Diff}(\ol \DD,p)
\;\cong\;
\mathrm{Diff}^{(n)}(\ol \DD,p),
\]
where \(\mathbb Z_n\) acts by postcomposition with the rotations \(r_{2\pi k/n}\). Endowing \(\mathrm{Diff}^{(n)}(\ol \DD,p)\) with the corresponding quotient smooth structure makes this identification a diffeomorphism and \(v_n\) an \(n\)-fold covering.

Furthermore, for \(\varphi_1,\varphi_2\in\mathrm{Diff}(\ol \DD,p)\),
\[
v_n(\varphi_1\circ\varphi_2)
=
v_n(\varphi_1)\circ\varphi_2,
\]
and, for \(r_\theta(z)=e^{i\theta}z\),
\[
v_n(r_\theta\circ\varphi)
=
r_{n\theta}\circ v_n(\varphi).
\]
\end{proposition}
\begin{proof}
    Surjectivity follows from the definition of \(\mathrm{Diff}^{(n)}(\ol \DD,p)\). For injectivity, suppose that
\[
v_n(\varphi_1)=v_n(\varphi_2).
\]
For \(z\ne p\), this implies
\[
\left(\frac{\varphi_2(z)}{\varphi_1(z)}\right)^n=1.
\]
The ratio is continuous on the connected punctured disk and takes values in the discrete set of \(n\)-th roots of unity. It is therefore constant, so
\[
\varphi_2=r_{2\pi k/n}\circ\varphi_1
\]
for some \(k\in\{0,\ldots,n-1\}\). This equality also holds at \(p\), since both diffeomorphisms fix \(p\). Conversely,
\[
u_n\circ r_{2\pi k/n}=u_n,
\]
so such rotations preserve the fibres of \(v_n\). Thus the fibres are exactly the \(\mathbb Z_n\)-orbits, proving the asserted bijection. 

The remaining identities follow directly from
\[
v_n(\varphi_1\circ\varphi_2)
=
u_n\circ\varphi_1\circ\varphi_2
=
v_n(\varphi_1)\circ\varphi_2
\]
and
\[
v_n(r_\theta\circ\varphi)
=
u_n\circ r_\theta\circ\varphi
=
r_{n\theta}\circ u_n\circ\varphi
=
r_{n\theta}\circ v_n(\varphi).
\]

\end{proof}

\begin{proposition}
    We have
    $$
{\rm Hyp}^{(n)}(\ol{\DD},p) \cong S^1 \backslash {\rm Diff}^{(n)}(\ol{\DD}, p), \hskip 0.3cm
{\rm Teich}^{(n)}(\ol{\DD}, p) \cong
S^1 \backslash {\rm Diff}^{(n)}(\ol{\DD}, p)/{\rm Diff}(\ol{\DD}; \partial \ol{\DD}, p).
    $$
\end{proposition}

\begin{proof}
    For the first statement, we compute
    $$
{\rm Hyp}^{(n)}(\ol{\DD}, p) \cong S^1 \backslash {\rm Diff}(\ol{\DD}, p) \cong S^1 \backslash {\rm Diff}^{(n)}(\ol{\DD}, p),
    $$
    where we have used the fact that $v_n$ induces a diffeomorphism between the quotients by rotations. For the second statement, we have
    $$
{\rm Teich}^{(n)}(\ol{\DD}, p) = {\rm Hyp}^{(n)}(\ol{\DD}, p)/{\rm Diff}(\ol{\DD}; \partial \ol{\DD}, p) =
S^1 \backslash {\rm Diff}^{(n)}(\ol{\DD}, p)/{\rm Diff}(\ol{\DD}; \partial \ol{\DD}, p)\cong S^1\backslash{\rm Diff}^{(n)}(S^1).
    $$
    Hence, by Proposition \ref{circle-n-folds} we obtain
    \begin{equation}
        {\rm Teich}^{(n)}(\ol{\DD},p)\cong S^1\backslash{\rm Diff}^{}(S^1). 
    \end{equation}
\end{proof}
\begin{rmk}
Note that by the preceding construction, the coadjoint orbit admits a geometric interpretation in terms of ${\rm Diff}(\ol \DD,p)$. We have the following identification:
\begin{equation}
\mathcal{O}^{(n)}\cong {\rm PSU(1,1)}\backslash{}{\rm Diff}^{(n)}(\ol{\DD},p)/{\rm Diff}^{(n)}(\ol{\DD},\partial \ol{\DD},p).
    \end{equation}
\end{rmk}

\section{Symplectic structures}
In this Section, we describe natural 2-forms on Teichmüller spaces ${\rm Teich}^\beta(\ol{\DD})$, and discuss their relation to Virasoro coadjoint orbits.

\subsection{Coframe formalism and the Atiyah-Bott 2-form}
We briefly recall the coframe formalism for hyperbolic metrics \cite{AM} and  the equivalent language of flat connections, and describe the definition and basic properties of  Atiyah-Bott  
2-forms \cite{atiyah-bott}.

\subsubsection{Cartan formalism and flat connections}\label{cartan}
Let $\Sigma$ be an oriented 2-manifold with a Riemannian metric $g$. On a chart $U \subset \Sigma$ one can choose an orthonormal coframe $(\alpha_1, \alpha_2)$ such that
$$
g = \alpha_1^2 + \alpha_2^2.
$$
Then, Cartan's moving frame formalism implies that there is a unique 1-form
$\sigma \in \Omega^1(U)$ (the spin connection) verifying the structure equations
$$
d \alpha_1 = \sigma \wedge \alpha_2, \hskip 0.3cm
d\alpha_2 = - \sigma \wedge \alpha_1.
$$
Recall that the metric $g$ is hyperbolic if and only if the ${\rm sl}(2, \mathbb{R})$-connection
$$
A = \frac{1}{2}
\left(
\begin{array}{cc}
\alpha_2 & \alpha_1 + \sigma \\
\alpha_1 - \sigma & - \alpha_2
\end{array}
\right)
$$
is flat. That is, 
$$
F_A= dA + A^2 =0.
$$
Let $H = {\rm SO}(2)$ represented by matrices
$$
h_\theta = 
\left(
\begin{array}{cc}
\cos(\theta/2) &  -\sin(\theta/2) \\
 \sin(\theta/2) & \cos(\theta/2)
\end{array}
\right).
$$
Different coframes over $U$ are related by coframe rotations
$\theta: U \to {\rm SO}(2)$ such that
$$
\alpha_1^\theta = \alpha_1 \cos(\theta) - \alpha_2 \sin(\theta), 
\hskip 0.3cm
\alpha_2^\theta = \alpha_1 \sin(\theta) + \alpha_2 \cos(\theta).
$$
This corresponds to a gauge transformation of the connection $A$:
\begin{equation}     \label{eq:gauge}
A^\theta = h_{\theta}^{-1} A h_{\theta} + h_{\theta}^{-1} dh_{\theta},
\end{equation}
and of the spin connection
$$
\sigma^\theta = \sigma - d\theta.
$$

Choose an invariant scalar product $(\cdot, \cdot)$ on the Lie algebra ${\rm sl}(2, \mathbb{R})$, for instance
$$
(a, b) = {\rm Tr}(ab).
$$
Let $P\to \Sigma$ be a principal ${\rm PSL}(2, \mathbb{R})$ bundle over $\Sigma$. In local charts, flat ${\rm sl}(2, \mathbb{R})$ connections on $\Sigma$ are represented by 1-forms $A\in \Omega^1(U, {\rm sl}(2, \mathbb{R}))$. Since ${\rm PSL}(2, \mathbb{R})$ is homotopically equivalent to ${\rm SO}(2)$, one can choose transformations between connection 1-forms as in \eqref{eq:gauge}.

The space of flat connections on $P$ carries a natural Atiyah-Bott 2-form \cite{atiyah-bott}
\begin{equation} \label{eq:AB}
\omega_{\rm AB} = \frac{1}{2} \int_\Sigma (\delta A, \delta A),
\end{equation}
where $\delta$ stands for the de Rham differential on the space of connections. While the 1-forms $A$ are only defined in local charts, it is easy to check that the right hand side of \eqref{eq:AB} is well defined under the gluing of charts.

For surfaces without boundary, the Atiyah-Bott 2-form \eqref{eq:AB} descends to the corresponding Teichmüller space, and induces a symplectic (Weil-Petersson \cite{wp_form}) structure.

\subsubsection{Surfaces with boundary}
When the surface $\Sigma$ has a boundary $\partial \Sigma \cong S^1$, and the flat connection $A$ is induced by a choice of a coframe for a hyperbolic $0$-metric, the integral in \eqref{eq:AB} may diverge at $\partial \Sigma$. In order to deal with this issue, introduce a local coordinate system $(\varrho, \phi)$ near the boundary. Here $\varrho$ is a local boundary defining function and $\phi$ is a coordinate on the boundary. Recall that a coframe $(\alpha_1, \alpha_2)$ is called adapted if
$$
\alpha_2 = \frac{d \varrho}{\varrho} + O(\varrho^0).
$$
Informally, one can imagine that in the dual frame $(e_1, e_2)$ the vector field $e_1$ is tangent to the boundary and $e_2$ is orthogonal to the boundary.

Denote by $\mathcal{A}^{\rm adp}$ the space of flat connections locally induced by adapted coframes of hyperbolic $0$-metrics. It turns out (see \cite{AM}) that when restricted to $\mathcal{A}^{\rm adp}$, the integral on the right hand side of \eqref{eq:AB} converges, and equation \eqref{eq:AB} defines a closed 2-form on $\mathcal{A}^{\rm adp}$. However, in contrast to the case of closed surfaces, $\omega_{\rm AB}$ does not descend to the Teichmüller space ${\rm Teich}(\Sigma)$.

Recall that the 1-form $A$ near the boundary $\partial \Sigma$ is in the Drinfeld-Sokolov \cite{drinfeld-sokolov} form if it reads
\begin{equation}     \label{eq:DS}
A = \left(
\begin{array}{cc}
\frac{1}{2} \frac{d\varrho}{\varrho} + O(\varrho^0) & \frac{d\phi}{\varrho} + O(\varrho^0) \\
-\varrho T(\phi) d\phi + O(\varrho^2) & -\frac{1}{2} \frac{d\varrho}{\varrho} + O(\varrho^0)
\end{array}
\right).
\end{equation}
We denote the space of flat connections in the Drinfeld-Sokolov form by 
$\mathcal{A}^{\rm DS}$ \footnote{More precisely, let $B$ be a singular gauge tranformation of $A$, that is
\[
B=r^{-1}Ar+r^{-1}dr,
\qquad
r=\begin{pmatrix}
\varrho^{-1/2}&0\\
0&\varrho^{1/2}
\end{pmatrix}.
\]
Then \(A\) is in Drinfeld–Sokolov form if \(B\) extends smoothly to the boundary and its pullback under the inclusion \(\iota:\partial\Sigma\hookrightarrow\Sigma\) satisfies
\[
\iota^*B=
\begin{pmatrix}
0&1\\
-T(\varphi)&0
\end{pmatrix}d\varphi.
\]
}.

The following theorem summarizes several results of \cite{AM}:
\begin{theorem}\label{AMtheorem}
(i) For each class of hyperbolic $0$-metrics $[g] \in {\rm Teich}(\Sigma)$ one can choose a metric $g$ and an adapted coframe $(\alpha_1, \alpha_2)$ such that the corresponding 1-form $A$ is in the Drinfeld-Sokolov form. 

(ii) The Atiyah-Bott form \eqref{eq:AB} descends from 
    $\mathcal{A}^{\rm DS}$ to a symplectic form on ${\rm Teich}(\Sigma)$. 

    (iii) The action of $\widetilde{\rm Diff}(S^1)$ (or of ${\rm Diff}(S^1)$) on ${\rm Teich}(\Sigma)$ is a Virasoro action with moment map
    $$
\mu([g]) = \partial_\phi^2 + T(\phi),
    $$ 
    where $T(\phi)$ is the coefficient in \eqref{eq:DS}.
\end{theorem}

As an example, consider the Poincaré metric $g_\DD$ on the disk, and the coframe
$$
\alpha_1 = \frac{2 \rho d\phi}{1-\rho^2}, \hskip 0.3cm
\alpha_2 =- \frac{2 d\rho}{1-\rho^2}.
$$
The corresponding spin connection is
$$
\sigma =  \frac{1+ \rho^2}{1- \rho^2} \, d\phi,
$$
and the connection $A$ is given by
$$
A = 
\left(
\begin{array}{cc}
- \frac{ d\rho}{1-\rho^2}    &  \frac{1}{2} \frac{1 + \rho}{1-\rho} \, d\phi \\
-\frac{1}{2} \frac{1-\rho}{1+\rho} \, d\phi     & \frac{d\rho}{1-\rho^2}
\end{array}
\right).
$$
By using the boundary defining function $\varrho = 1-\rho$, it is easy to check that this connection is in the Drinfeld-Sokolov form with $T(\phi)= \frac{1}{4}$.

\subsection{Teichmüller spaces ${\rm Teich}^\beta(\ol{\DD}, p)$}
In this Section, we construct a natural closed \(2\)-form on \(\operatorname{Teich}^{\beta}(\overline{\mathbb D},p)\) and establish its relation to Virasoro coadjoint orbits. We first describe the corresponding coframes and flat connections, and compute the two boundary contributions to the Atiyah–Bott form. We then impose the Drinfeld–Sokolov condition at the ideal boundary and show that, after subtracting the contribution at the marked point, the resulting form descends to Teichmüller space.

\subsubsection{Coframes and flat connections on $\ol\DD \backslash p$}

Consider the hyperbolic metric 
$$
g_\DD^\beta = \frac{4 \beta^2 |z|^{2 \beta -2} |dz|^2}{(1- |z|^{2\beta})^2},
$$
the corresponding coframe $\alpha^\beta$
$$
\alpha^\beta_1 = \frac{2 \beta \rho^{\beta}}{(1 - \rho^{2\beta})}\, d\phi, \hskip 0.3cm
\alpha^\beta_2 = -\frac{2 \beta \rho^{\beta-1 } d\rho}{(1 - \rho^{2\beta})},
$$
the spin connection
$$
\sigma =  \beta \frac{1 + \rho^{2\beta}}{1-\rho^{2\beta}} \, d\phi,
$$
and the flat ${\rm sl}(2, \mathbb{R})$ connection
$$
A^\beta = 
\left(
\begin{array}{cc}
 -\beta\frac{  \rho^{\beta-1}d\rho}{1-\rho^{2\beta}}    &  \frac{\beta}{2} \frac{1 + \rho^\beta}{1-\rho^\beta} \, d\phi \\
-\frac{\beta}{2} \frac{1-\rho^\beta}{1+\rho^\beta} \, d\phi     & \beta\frac{  \rho^{\beta-1}d\rho}{1-\rho^{2\beta}}
\end{array}
\right).
$$

\begin{proposition}
    The connection $A^\beta$ is in the DS form with $T^\beta(\phi) = \frac{1}{4} \beta^2$.
\end{proposition}

\begin{proof}
  With the boundary defining function 
    $\varrho \sim 1- \rho$, we obtain
    $$
-\beta\frac{ \rho^{\beta-1} d\rho}{1-\rho^{2\beta}}  = \frac{1}{2} \, \frac{d \varrho}{\varrho} + O(\varrho^0),
    $$
and

    $$
\frac{\beta}{2} \frac{1 + \rho^\beta}{1-\rho^\beta} \, d\phi =
\frac{d \phi}{\varrho} + O(\varrho^0), \hskip 0.3cm
-\frac{\beta}{2} \frac{1-\rho^\beta}{1+\rho^\beta} \, d\phi  = -\frac{\beta^2}{4} \varrho + O(\varrho^2).
    $$ Therefore the connection $A^\beta$ acquires the DS form
    \begin{equation}\label{Ads}
 A^\beta = \left(
\begin{array}{cc}
\frac{1}{2} \frac{d\varrho}{\varrho}+O(\varrho^0) & \frac{d\phi}{\varrho}+O(\varrho^0)  \\
-\varrho \frac{\beta^2}{4}d\phi+O(\varrho^2) & -\frac{1}{2} \frac{d\varrho}{\varrho} +O(\varrho^0)
\end{array}
\right). 
    \end{equation}
    \end{proof}

Define the group
$$
\mathcal{G} = \{ (\varphi, \theta) \in {\rm Diff}(\ol{\DD};  p) \ltimes {\rm Map}(\ol{\DD}, {\rm SO}(2));
\theta|_{\partial\ol{\DD}}=0\}
$$
as a semidirect product of the group of orientation preserving diffeomorphisms of $\ol{\DD}$ preserving the boundary and the point $p$, and the gauge group of coframe rotations trivial at the boundary $\partial\ol{\DD}$. This group acts on coframes and we denote the orbit of $\alpha^\beta$ under its action by
$$
{\rm Cof}^\beta(\DD) = \mathcal{G} \cdot \alpha^\beta.
$$
Fixing a standard frame at $p$, the differential $d_p\varphi$ can be represented as a matrix
\begin{equation}\label{dphi}
    L=d_p \varphi =
\left(
\begin{array}{cc}
   a  &  b \\
   c  &  d
\end{array}
\right) \in {\rm GL^+}(2, \mathbb{R}), \qquad  \Delta=ad-bc>0
\end{equation}   
Let $\nabla \subset \rm GL^+(2,\mathbb{R})$ be the subgroup of upper triangular matrices with positive diagonal entries. We define the subgroup of $\mathcal{G}$
\begin{equation*}
    \mathcal{G}_0=\{ (\varphi, \theta) \in \mathcal{G}: d_p\phi \in \nabla\}.
\end{equation*}
 
\begin{proposition}\label{coframe stabilizers}
    Stabilizers of the natural action of $\mathcal{G}$ on ${\rm Cof}^\beta(\DD)$ are isomorphic to ${\rm SO}(2)$, all elements $\alpha \in {\rm Cof}^\beta(\DD)$ are adapted coframes, the natural map 
    ${\rm Cof}^\beta(\DD) \to {\rm Hyp}^\beta(\DD; p)$ given by
    $$
\alpha=(\alpha_1, \alpha_2) \mapsto g=\alpha_1^2 + \alpha_2^2
    $$
    is surjective. 
\end{proposition}

\begin{proof}
    First, let $g \in {\rm Hyp}^\beta(\DD, p)$. There is a diffeomorphism
    $\varphi \in {\rm Diff}(\ol{\DD};  p)$ such that 
    $g=\varphi^* g^\beta_{\DD}$, where $g^\beta_{\DD}$ is the standard metric. Then, taking the image $\alpha= \varphi^* \alpha^\beta$ of the standard coframe $\alpha^\beta$ one obtains a coframe corresponding to $g$. Furthermore, since coframe rotations preserve metrics, ${\rm Hyp}^\beta(\DD)$ is the full image of ${\rm Cof}^\beta(\DD)$.

    Next, recall that under the action of diffeomorphisms the stabilizer of $g^\beta_{\DD}$ is equal to 
    $S^1 \cong {\rm SO}(2) \subset {\rm Diff}(\ol{\DD}; p)$. By inspection, this subgroup also preserves the coframe $\alpha^\beta$. Since the action of coframe rotations on coframes is free, the stabilizer of $\alpha^\beta$ under the action of $\mathcal{G}$ is also equal to ${\rm SO}(2)$, and stabilizers of all other coframes $\alpha \in {\rm Cof}^\beta(\DD)$ are isomorphic to ${\rm SO}(2)$.

    Finally, we compute
    $$
\varphi^* \frac{d \varrho}{\varrho} = \frac{d (\lambda(\phi) \varrho + O(\varrho^2))}{\lambda(\phi) \varrho + O(\varrho^2)} =
\frac{d \varrho}{\varrho} + O(\varrho^0).
    $$
    Since for elements of $\mathcal{G}$ we have
    $\theta|_{\partial\ol{\DD}}=0$, the property above is preserved by coframe rotations, and all coframes in ${\rm Cof}^\beta(\DD)$ are adapted.
\end{proof}

We denote the set of flat ${\rm sl}(2, \mathbb{R})$ connections on $\DD\backslash \{ p\}$ corresponding to coframes in ${\rm Cof}^\beta(\DD)$ by $\mathcal{A}^\beta$. Proposition \ref{coframe stabilizers} allows to write $\mathcal{A}^\beta$ as the homogenous space $\mathcal{A}^\beta\cong \S^1\backslash\mathcal{G}$. The $\mathcal{G}$-action on the coframes induces  the orbit map 
 \begin{align*}
    q: \mathcal{G}&\rightarrow \mathcal{A}^\beta \\
        (\varphi, \theta) &\mapsto h_\theta^{-1}\varphi^*A^\beta h_\theta + h_\theta^{-1}d h_\theta.
    \end{align*}
Each fiber  of $q$ is an $S^1$ orbit.

\begin{proposition}\label{global_section}
For every $A\in\mathcal A^\beta$, there exists a unique
pair $(\varphi_A,\theta_A)\in\mathcal G_0$ such that
\[
q(\varphi_A,\theta_A)=A.
\]
The resulting map
\[
\sigma:\mathcal A^\beta\longrightarrow\mathcal G_0,
\qquad A\longmapsto(\varphi_A,\theta_A),
\]
is a smooth global section of $q$.
\end{proposition}

\begin{proof}
Let $A\in\mathcal A^\beta$ and choose
$(\varphi,\theta)\in q^{-1}(A)$.
Consider the QR decomposition of a $\rm GL^+(2,\mathbb{R})$ element
\[
d_p\varphi=QT,\qquad
Q\in\mathrm{SO}(2),\quad T\in\nabla.
\]
Since rigid rotations preserve $A^\beta$, the pair
$(r_Q^{-1}\circ\varphi,\theta)$ belongs to the same
fiber as $(\varphi, \theta)$. Furthermore,
\[
d_p(r_Q^{-1}\circ\varphi)
=Q^{-1}d_p\varphi=T\in\nabla,
\]
so this pair belongs to $\mathcal G_0$.
Uniqueness of this pair follows from the uniqueness of the QR decomposition. Hence, the restriction $q|_{\mathcal{G}_0}$ is a bijection. Its inverse provides us with the smooth  global section
\begin{equation*}
\sigma(A)=q|_{\mathcal{G}_0}^{-1}: \mathcal{A}^\beta \rightarrow \mathcal{G}_0 \subset \mathcal{G}, \qquad \sigma=(\phi_A,\theta_A),
\end{equation*}
and $q \circ \sigma=\rm id_{\mathcal{A}^\beta}.$
\end{proof}

\begin{proposition}\label{pure_gauge}
    All connections $A \in \mathcal{A}^\beta$ admit a form
    $A = G^{-1} dG$ for a quasi-periodic function
    $$
G(\rho, \phi+2\pi) = M_\beta G(\rho, \phi),
    $$
    with the monodromy matrix $M_\beta = h_{-2 \pi \beta}$.
\end{proposition}

\begin{proof}
    By inspection, the connection $A^\beta$ satisfies the required property with 
    $$
G^\beta=
\left(
\begin{array}{cc}
\cos\left(\frac{\beta}{2}\phi\right) & \sin\left(\frac{\beta}{2}\phi\right) \\
-\sin\left(\frac{\beta}{2}\phi\right) & \cos\left(\frac{\beta}{2}\phi\right)
\end{array}
\right)
\left(
\begin{array}{cc}
\left(\frac{1-\rho^\beta}{1+\rho^\beta}\right)^{\frac{1}{2}} & 0 \\
0 & \left(\frac{1+\rho^\beta}{1-\rho^\beta}\right)^{\frac{1}{2}}
\end{array}
\right).
    $$
    Applying an element $(\varphi, \theta) \in \mathcal{G}$, we obtain
    $A = G^{-1}dG$ with $G= \varphi^*(G^\beta)h_\theta$. From this explicit formula, it is clear that $G$ satisfies the same quasi-periodicity condition as  $G^\beta$.
\end{proof}

For future use, we will introduce a real blow-up $\DD_p$ of $\DD$ at the point $p$. For a diffeomorphism $\varphi$, denote $\tilde{\phi}=\varphi^* \phi$. Then, restricting to $\rho=0$ and putting $h_0=h(\rho=0)$ we obtain
$$
A(\rho=0) = 
h_0^{-1} (\varphi^* A^{\beta}(\rho=0))h_0 =
\left(
\begin{array}{cc}
0 & \frac{\beta}{2} d \tilde{\phi} \\
-\frac{\beta}{2} d \tilde{\phi} & 0
\end{array}
\right).
$$
Here we have taken into account that the conjugation by $h_0$ doesn't change this connection.

\subsubsection{Computation of the symplectic form}
 The section $\sigma$ in Proposition \ref{global_section} identifies the total space with $\mathcal{A}^\beta \times S^1$. The form $\delta GG^{-1}$ is a one-form on $\mathcal{A}^\beta \times S^1$ with values in quasi-periodic functions on $\DD_p$, with the periodicity condition \footnote{$\delta$ is the de Rham differential on the total space.} ,
$$
\delta G G^{-1}(\rho, \phi+2\pi) = M_\beta 
(\delta G G^{-1})(\rho, \phi) M_\beta^{-1}.
$$
On the other hand, we have $G^{-1}dG=A$, and
\begin{equation*}
    q^*(\delta A)=\delta(G^{-1}dG)=G^{-1}d(\delta G G^{-1})G.
\end{equation*}
From now on, we work along the section \(\sigma\) of Proposition \ref{global_section}. In particular, expressions involving \(\delta GG^{-1}\) are understood after pullback along \(\sigma\), which we leave implicit in the notation.

For $0 < \epsilon < \rho < 1$, let 
$$
\DD_p(\epsilon, \rho) = [\epsilon, \rho] \times S^1 \subset \DD_p,
$$
and define the Atiyah-Bott form on $\mathcal{A}^\beta$:
$$
\omega_{\rm AB}^{\epsilon, \rho}  =  \frac{1}{2} \int_{\DD_p(\epsilon, \rho)} {\rm Tr}  \, \delta A \wedge \delta A.
$$
We put $G_{\rho}=G(\rho, \phi)$. Note that $G_0(\phi)=G(0, \phi)$ is well defined. 

\begin{proposition}\label{w-stokes}
   For $0<\epsilon<\rho<1$, the Atiyah-Bott form reads

   $$
\begin{array}{lll}
\omega_{\rm AB}^{\epsilon, \rho} & = & \frac{1}{2} \int_{S^1}
{\rm Tr} (\delta G_{\rho} G_{\rho}^{-1}) d(\delta G_{\rho} G_{\rho}^{-1})\\
& -& 
\frac{1}{2} \int_{S^1}
{\rm Tr} (\delta G_\epsilon G_\epsilon^{-1}) d(\delta G_\epsilon G_\epsilon^{-1}).
\end{array}$$
\end{proposition}

\begin{proof}
    First, we note that
    $$
    \begin{array}{lll}
   \omega_{\rm AB}^{\epsilon, \rho} & = & \frac{1}{2}\int_{\DD_p(\epsilon,\rho)}{\rm Tr} \delta A \wedge \delta A \\& = &
    \frac{1}{2} \int_{\DD_p(\epsilon,\rho)} {\rm Tr} \, (\delta(G^{-1}dG) \wedge \delta(G^{-1}dG)) \\
    & = &
    \frac{1}{2} \int_{\DD_p(\epsilon,\rho)} {\rm Tr} \, 
    d(\delta G G^{-1}) \wedge d (\delta G G^{-1}) \\
    & = &
    \frac{1}{2}  \int_{\DD_p(\epsilon,\rho)} d \, {\rm Tr} \, 
   (\delta G G^{-1}) \wedge d (\delta G G^{-1}) \\
    & = &
    \frac{1}{2} \int_{\partial \DD_p(\epsilon,\rho)}
    {\rm Tr} \, 
    (\delta G G^{-1}) \wedge d (\delta G G^{-1}) \\
    & = & \frac{1}{2} \int_{S^1}
{\rm Tr} (\delta G_{\rho} G_{\rho}^{-1}) d(\delta G_{\rho} G_{\rho}^{-1})\\
& -& 
\frac{1}{2} \int_{S^1}
{\rm Tr} (\delta G_\epsilon G_\epsilon^{-1}) d(\delta G_\epsilon G_\epsilon^{-1}),
   \end{array}
$$
where we applied Stokes’ theorem on the annulus \(\DD_p(\epsilon,\rho)\), to express the Atiyah–Bott form as the sum of two circle boundary contributions $C_\epsilon$ and  $C_\rho$ with opposite orientations (Figure \ref{stokes picture}). We also used the fact that despite $\delta G G^{-1}$ being quasi-periodic, the 1-form ${\rm Tr} \, 
    (\delta G G^{-1}) \wedge d (\delta G G^{-1})$ descends to the $\DD_p(\rho)$.
\end{proof}
\begin{center}
\begin{figure}[H]
\begin{tikzpicture}[>=stealth, scale=0.7, transform shape]
\tikzset{
  boundary/.style={very thick, black},
  contour/.style={
    very thick, red!75!black,
    dash pattern=on 5pt off 2.5pt
  },
  edgeR/.style={very thick, red!75!black},
  lbl/.style={font=\normalsize, text=black}
}

\draw[boundary] (0,0) circle (1.75);

\draw[contour] (0,0) circle (1.5);
\draw[contour] (0,0) circle (0.35);

\draw[edgeR,->]
  (260:1.5)
  arc[start angle=260, end angle=280, radius=1.5];

\draw[edgeR,->]
  (290:0.35)
  arc[start angle=290, end angle=250, radius=0.35];

\fill[black] (0,0) circle (1.2pt);
\node[lbl] at (0,0.65) {$C_\epsilon$};
\node[lbl] at (0,-1.05) {$C_\rho$};
\end{tikzpicture}
\caption{The annulus $\DD_p(\epsilon, \rho)$}\label{stokes picture}
\end{figure}
\end{center} 

In the neighborhood of $\partial \ol{\DD}$, the diffeomorphism 
$\varphi: (\rho, \phi) \mapsto (\tilde{\rho}, \tilde{\phi})$ and the rotation angle
$\theta$ take the following form in local coordinates $(\varrho, \phi)$:
$$
\tilde{\rho} = 1 - \varrho g(\phi) + O(\varrho^2), \hskip 0.3cm
\tilde{\phi} = f(\phi) + O(\varrho), \hskip 0.3cm
\theta=\varrho h(\phi),
$$
where $\varrho = 1- \rho$ is a boundary defining function, $g: S^1 \to \mathbb{R}_+, h: S^1 \to \mathbb{R}$ are two functions of $\phi$, and $f \in {\rm Diff}(S^1)$ is an orientation preserving diffeomorphism of the circle. 

\begin{theorem}\label{2form-tip-boundary}
    The family of 2-forms $\omega_{\rm AB}^{\epsilon, \rho}$ has a limit
    $\omega_{\rm AB} = \lim_{\substack{\epsilon\to0\\ \rho\to1}}\omega_{\rm AB}^{\epsilon, \rho}$ given by
    $$
\omega_{\rm AB} = \omega_{\rm in} + \omega_{\rm out},
$$
  where       $$
    \begin{array}{lll}
    \omega_{\rm in} & = & \frac{\pi \beta^2}{4(ad-bc)^2}\bigl(d\,\delta a-c\,\delta b+b\,\delta c-a\,\delta d\bigr)\wedge \bigl(-b \delta a+a \delta b + d \delta c - c\delta d \bigr),\\
\omega_{\rm out} & = &  \int_{S^1} \left( \frac{1}{4} \delta(\log g) \wedge \delta((\log g)') - \frac{1}{2}\delta\left( \frac{f'}{g}\right) \wedge \delta h -
\frac{\beta^2}{4} \delta f \wedge \delta(f')\right).
\end{array}
    $$
\end{theorem}

\begin{proof}
We proceed in three steps for the proof.

\vspace{3mm}
\paragraph{Step 1.}

Let us define the two contributions. Evaluating the expression in Proposition \ref{w-stokes} in the limit $\epsilon \mapsto 0$, the inner boundary contribution gives
\begin{align*}
\omega_{\rm in}
&:=-\lim_{\epsilon\rightarrow0}\frac12\int_{S^1}
  \text{Tr}\!\left[
  (\delta G_\epsilon G_\epsilon ^{-1})\wedge d(\delta G_\epsilon G_\epsilon^{-1})\right]\\
&=-\frac12\int_0^{2\pi}
  \text{Tr}\!\left[
  (\delta G_0G_0^{-1})\wedge
  \partial_\phi(\delta G_0G_0^{-1})\right]\,d\phi.
\end{align*}

We then define the outer boundary contribution by 
\[
\omega_{\rm out}
:=\lim_{\rho\rightarrow1}\frac12\int_{S^1}
\rm Tr\!\left[
  (\delta G_\rho G_\rho^{-1})\wedge d(\delta G_\rho G_\rho^{-1})\right],
\]
and we will see in Step 3 that this limit exists.
Therefore, 
\begin{equation*}
    \omega_{\rm AB}=\lim_{\substack{\epsilon\to0\\ \rho\to1}} \omega_{\rm AB}^{\epsilon, \rho}=\omega_{\rm in} + \omega_{\rm out}.
\end{equation*}

\vspace{3mm}

\paragraph{Step 2.}We compute the inner boundary contribution. 
Recall $L$ introduced in \eqref{dphi},
\[
L=d_p\varphi=\begin{pmatrix}a&b\\c&d\end{pmatrix}, \qquad \Delta=ad-bc>0.
\]
The induced angular action under such transformation is
$\widetilde\phi=\arg (L(\cos\phi,\sin\phi)^T)$, or equivalently \footnote{$(u_1,u_2)^T$ denotes the transpose of the vector $u=\begin{pmatrix}
    u_1 \\u_2
\end{pmatrix}$.},
\begin{equation*}
    \tan\widetilde\phi=\frac{c+d\tan\phi}{a+b\tan\phi}.
\end{equation*}

For $J=\begin{pmatrix}0&-1\\1&0\end{pmatrix}$ and $h_t=\exp(tJ/2)$,
we obtain
\[
G_0=h_{-\beta\widetilde\phi}h_{\theta_0},\qquad
\delta G_0G_0^{-1}
=\frac J2(-\beta\,\delta\widetilde\phi+\delta\theta_0),
\]
where $\theta_0$ is the coframe rotation at $p$, independent of $\phi$.
The inner boundary term is then
\begin{equation}\label{inner integral}
\begin{aligned}
\omega_{\rm in}
&=-\frac12\int_0^{2\pi}\Tr\!\left[
(\delta G_0G_0^{-1})\wedge\partial_\phi(\delta G_0G_0^{-1})
\right]\,d\phi\\
&=\frac{\beta^2}{4}\int_0^{2\pi}
\delta\widetilde\phi\wedge\partial_\phi(\delta\widetilde\phi)\,d\phi.
\end{aligned}
\end{equation}
The term involving $\delta\theta_0$
vanishes since $\int_0^{2\pi}\partial_\phi(\delta\widetilde\phi)\,d\phi=0$.

To evaluate \eqref{inner integral} in terms of the entries of $L$, we compute 
\begin{equation}
\label{eq:theta}
\delta L\,L^{-1}
=\frac1\Delta\begin{pmatrix}
d\,\delta a-c\,\delta b&-b\,\delta a+a\,\delta b\\
d\,\delta c-c\,\delta d&-b\,\delta c+a\,\delta d
\end{pmatrix}.
\end{equation}
Write $v=L(\cos\phi,\sin\phi)^T
=s(\cos\widetilde\phi,\sin\widetilde\phi)^T$, with $s>0$.
At fixed $\phi$, $\delta v=(\delta L)L^{-1}v$.
Its scalar product with the unit angular vector
$(-\sin\widetilde\phi,\cos\widetilde\phi)^T$ is
$s\,\delta\widetilde\phi$. Therefore to compute $\delta\widetilde{\phi}$ it suffices to do:
\begin{equation}\label{eq:angle-variation}
\begin{aligned}
\delta\widetilde\phi
&=(-\sin\widetilde\phi,\cos\widetilde\phi)\,
(\delta L L^{-1})\begin{pmatrix}\cos\widetilde\phi\\\sin\widetilde\phi\end{pmatrix}.\\
\end{aligned}
\end{equation}
 Integrating then the trigonometric terms in
\eqref{inner integral} gives
\begin{equation}\label{eq:tip-trace}
\omega_{\rm in}
=\frac{\pi\beta^2}{4}\Tr(J \delta L L^{-1}\wedge\delta L L^{-1}).
\end{equation}
In terms of the entries of $L$, this is
\begin{equation}\label{eq:tip-coordinates}
{\begin{aligned}
\omega_{\rm in}=\frac{\pi\beta^2}{4(ad-bc)^2}\,
&\bigl(d\,\delta a-c\,\delta b+b\,\delta c-a\,\delta d\bigr)\\[-1pt]
{}\wedge&\bigl(-b\,\delta a+a\,\delta b+d\,\delta c-c\,\delta d\bigr).
\end{aligned}}
\end{equation}

\vspace{3mm}
\paragraph{Step 3.}Finally, let us now compute the outer boundary contribution. 
Let $(\widetilde{\rho},\tilde{\phi}) $ be the image of a point $({\rho},{\phi}) \in \ol{\DD}$   under the  action of $\varphi$, and let $\theta$ denote the coframe rotations. Using the boundary defining function $\varrho=1-\rho$, we can write the expansions near the outer boundary:
\[
\widetilde\rho=1-\varrho g(\varphi) + O(\varrho^2),\qquad
\widetilde\varphi=f(\varphi) + O(\varrho),\qquad
\theta=\varrho h(\varphi) + O(\varrho^2).
\] 
Introduce the constant matrices
\[
K=\frac12\begin{pmatrix}0&-1\\1&0\end{pmatrix},\qquad
H=\begin{pmatrix}1&0\\0&-1\end{pmatrix},\qquad
S=\frac12\begin{pmatrix}0&1\\1&0\end{pmatrix}.
\]
Then we can write $G$, including the coframe rotations, as
\begin{equation}
G=e^{-\beta \tilde{\varphi} K}e^{-sH}e^{\theta K},
\end{equation}
with 
$s=\frac12\log\frac{1+\widetilde\rho^{\,\beta}}{1-\widetilde\rho^{\,\beta}}$.
 The radial variation is \footnote{The field-space differential $\delta$ acts on $f,g,h$, at fixed
$\varrho,\phi$.}
\begin{equation}\label{deltas}
\delta s
=
\frac{\beta\widetilde\rho^{\beta-1}}
     {1-\widetilde\rho^{2\beta}}\,
\delta\widetilde\rho.
\end{equation}
Writing $\xi=\delta G_\rho\,G_\rho^{-1}$, and using the identity
$e^{sH}Ke^{-sH}=K\cosh(2s)-S\sinh(2s)$, we get
\begin{equation*}
\begin{aligned}
\xi&=e^{-\beta \tilde{\varphi}K}Qe^{\beta \tilde{\varphi}K},\\
Q&=-\beta K\,\delta \tilde{\varphi}-H\,\delta s
 +\bigl(K\cosh(2s)+S\sinh(2s)\bigr)\delta\theta.
\end{aligned}
\end{equation*}
Differentiating yields (primes denote $\partial_\varphi$):
\[
\xi'=e^{-\beta \tilde{\varphi}K}\bigl(Q'-[\beta \tilde{\varphi}'K,Q]\bigr)e^{\beta \tilde{\varphi}K},
\qquad
\Tr(\xi\wedge\xi')=\Tr\!\left(Q\wedge\bigl(Q'-[\beta \tilde{\varphi}'K,Q]\bigr)\right).
\]
Putting $c=\cosh(2s)$, one obtains
\begin{equation}
\begin{aligned}
\Tr(\xi\wedge\xi')={}&
2\,\delta s\wedge\delta s'-\frac{\beta^2}{2}\,\delta \tilde{\varphi}\wedge\delta \tilde{\varphi}'
-\frac{1}{2}\,\delta\theta\wedge\delta\theta'\\
&+\frac{\beta c}{2}\bigl(\delta \tilde{\varphi}\wedge\delta\theta'
+\delta\theta\wedge\delta \tilde{\varphi}'\bigr)
+\frac{\beta c'}{2}\,\delta \tilde{\varphi}\wedge\delta\theta
-\beta \tilde{\varphi}'\,\delta c\wedge\delta\theta,
\end{aligned}
\end{equation}
where using integration by parts the terms involving $c'$ cancel. The remaining mixed terms combine into
$-\frac{1}{2}\beta \,\delta(c\tilde{\varphi}')\wedge\delta\theta$. Hence, 
\begin{equation}
\int_{S^1}\frac{1}{2}\Tr(\xi\wedge\xi')\,d\varphi\\
=\int_{S^1}\left[
\,\delta s\wedge\delta s'
-\frac{\beta^2}{4}\,\delta \tilde{\varphi}\wedge\delta \tilde{\varphi}'
-\,\frac{1}{4}\delta\theta\wedge\delta\theta'
-\frac{1}{2}\beta\,\delta(c\tilde{\varphi}')\wedge\delta\theta
\right]d\varphi.
\end{equation}
We now take the limit \(\varrho\to0\). The expansion
\(\widetilde\rho=1-\varrho g+O(\varrho^2)\) gives
\[
\delta\widetilde\rho=-\varrho\,\delta g+O(\varrho^2),
\qquad
\frac{\beta\widetilde\rho^{\beta-1}}
     {1-\widetilde\rho^{2\beta}}
=\frac{1}{2\varrho g}+O(1).
\]
Substituting into the radial variation \eqref{deltas}, we obtain
\[
\delta s
=
\left(\frac{1}{2\varrho g}+O(1)\right)
\left(-\varrho\,\delta g+O(\varrho^2)\right)
=
-\frac12\,\delta\log g+O(\varrho).
\]
Thus \(\delta s\) has a finite boundary limit, even though \(s\) itself diverges. Its contribution to the 2-form is
\[\lim_{\varrho \to 0}
\delta s\wedge\delta s'
=
\frac14\delta\log g\wedge(\delta\log g)'.
\]
 Note also the expansion
$
c=\cosh{(2s)}=\frac{1+\widetilde\rho^{\,2\beta}}{1-\widetilde\rho^{\,2\beta}}
 =\frac{1}{\beta\varrho g}+O(1).
$
Together with $\delta\theta=\varrho\delta h$ and $\delta \tilde{\varphi}=\delta f + O(\varrho)$ we obtain in the boundary limit:
\begin{align*}
\lim_{\varrho \to 0}(\beta\,\delta(cf')\wedge\delta\theta)
&= \,\delta\!\left(\frac{f'}g\right)\wedge\delta h,\qquad
\lim_{\varrho \to 0}\delta\theta\wedge\delta\theta'
=0.
\end{align*}
Finally, 
\begin{equation}\label{w1}
{\displaystyle
\omega_{\rm out}=\lim_{\varrho \to 0}\int_{S^1}\frac{1}{2}\Tr(\xi\wedge\xi')\,d\varphi\\=\int_{S^1}\left[
-\frac{\beta^2}{4}\,\delta f\wedge\delta f'
+\frac14\,\delta\log g\wedge(\delta\log g)'
-\frac{1}{2}\,\delta\!\left(\frac{f'}g\right)\wedge\delta h
\right]d\varphi,}
\end{equation}
as required.
\end{proof}

Note that \(\omega_{\rm in}\) and \(\omega_{\rm out}\) are closed \(2\)-forms on \(\rm{Cof}^{\beta}(\overline{\mathbb D})\). The form \(\omega_{\rm in}\) descends to the space of metrics, where it is the pullback of a symplectic form on the space of complex structures on \(T_p\mathbb D\), see Remark \ref{cpx_structures} below. The form \(\omega_{\rm out}\) does not descend in general.  

We now describe the  describe the geometric meaning of \(\omega_{\rm in}\), namely we will show that \(\omega_{\rm in}\) is proportional to the pullback of the Kirillov–Kostant–Souriau (KKS) symplectic form on the coadjoint orbit.

For a group $\mathcal{F}$, let $\xi \in \mathfrak{f}^*$ and denote by
$$
O_\xi =\{ {\rm Ad}^*_u \xi; u \in \mathcal{F}\}
$$
the corresponding coadjoint orbit. Denote by $\pi_\xi: u \mapsto {\rm Ad}^*_u \xi$ the canonical projection. Recall that $O_\xi$ carries a canonical KKS symplectic form $\omega_\xi$ which is characterized by the following property:
$$
\pi_\xi^* \omega_\xi = \frac{1}{2} \langle \xi, [u^{-1}\delta u, u^{-1} \delta u]\rangle.
$$

\begin{proposition}
   Let \(\mathcal O_\beta\) be the coadjoint orbit of
\[
\xi_\beta=\beta J,\qquad
J=\begin{pmatrix}0&-1\\1&0\end{pmatrix},
\]
under \(\mathrm{GL}^+(2,\mathbb R)\), using the trace pairing to identify the Lie algebra with its dual. Denote its KKS form by $\omega_{\mathrm{\xi_\beta}}$. Define
\[
\pi_\beta:\mathcal A^\beta\longrightarrow\mathcal O_\beta,
\qquad
\pi_\beta(A)=L(A)^{-1}\xi_\beta L(A),
\]
where \(L(A)=d_p\varphi_A\) is determined by the section of Proposition \ref{global_section}. Then
\[
\omega_{\rm in}=\frac{\pi\beta}{4}\,
\pi_\beta^*\omega_{\mathrm{\xi_\beta}.}
\]
\end{proposition}

\begin{proof}
Since $\pi_{\beta}(A)=\pi_{\xi_\beta}(L^{-1}(A))$, we apply the KKS pullback identity with \(u=L^{-1}\).
Then
\[
u^{-1}\delta u=L\,\delta(L^{-1})=-\delta L\,L^{-1}.
\]
Using $[\delta L\,L^{-1},\delta L\,L^{-1}]=2\delta L\,L^{-1}\wedge\delta L\,L^{-1}$ for
matrix-valued one-forms, we obtain
\[
\begin{aligned}
\frac{\pi\beta}{4}\pi_\beta^*\omega_{\mathrm{\xi_\beta}.}
&=\frac{\pi\beta}{4}\frac12\operatorname{Tr}\!\left(\xi_\beta[-\delta L\,L^{-1},-\delta L\,L^{-1}]\right)\\
&=\frac{\pi\beta^2}{4}\operatorname{Tr}(J\delta L\,L^{-1}\wedge\delta L\,L^{-1})
=\omega_{\rm in},
\end{aligned}
\]
as in \eqref{eq:tip-trace}.
\end{proof}

\begin{rmk}\label{cpx_structures}
The inner boundary contribution $\omega_{\rm in}$ has a natural interpretation
in terms of the complex structure at the marked point \eqref{eq:tip-trace}.
 Denote by $\mathcal J^+(T_p\mathbb{D})$ the space of linear complex structures
on $T_p\mathbb{D}$, with $J=\begin{pmatrix}0&-1\\1&0\end{pmatrix}$ being the standard one. For $g=\varphi^*g_D^\beta$,
the induced complex structure extends across $p$, and in the fixed
oriented frame we have
\[
    I_g(p)=L^{-1}JL,\qquad L=d_p\varphi,
    \qquad .
\]
By the preceding proposition, $\omega_{\rm in}$ is a constant multiple of
the pullback of the KKS form on $\mathcal J^+(T_p\mathbb{D})$ along
$\alpha\mapsto I_g(p)$, where $g$ is the metric defined by the coframe $\alpha$.
It therefore descends to the space of metrics.

Moreover, this description is compatible with the section defined in Proposition \ref{global_section}. Indeed, for the QR decomposition $L=QT$, with
$Q\in\mathrm{SO}(2)$ and $T\in\nabla$, we have
\[
    L^{-1}JL=T^{-1}JT,
\]
because $Q$ commutes with $J$. Two matrices
$T,T'\in\nabla$ determine the same complex structure if and only if
$T'=rT$ for some $r>0$. Consequently,
\[
    \nabla/\mathbb R_{>0}\longrightarrow\mathcal J^+(T_p\mathbb{D}),
    \qquad [T]\longmapsto T^{-1}JT,
\]
is a diffeomorphism, where $\mathbb R_{>0}$ denotes the subgroup of
positive scalar matrices. It follows then quite easily that the space $\mathcal J^+(T_p\mathbb{D})$ is also a coadjoint orbit of $\nabla$.
\end{rmk}

\subsubsection{The Drinfeld-Sokolov slice}

  For the form $\omega_{\rm out}$, recall that all connections in $\mathcal{A}^\beta$ can be brought to the Drinfeld-Sokolov (DS) slice \eqref{Ads}. More precisely, the DS form \cite{drinfeld-sokolov} is such that its boundary restriction has upper-right coefficient \(1\) and vanishing diagonal coefficients.
  This can be obtained by applying a singular gauge transformation to $A^\beta$:
\[
B\coloneqq r^{-1}A^\beta r+r^{-1}dr,
\qquad
r=
\begin{pmatrix}
\varrho^{-1/2}&0\\
0&\varrho^{1/2}
\end{pmatrix}.
\]
Then \(B\) extends smoothly to the boundary, and its boundary restriction $\iota^*B$ has upper-right coefficient \(1\) and vanishing diagonal coefficients, where \(\iota:\partial\ol{\DD}\hookrightarrow\DD\).

For the expansions
\[
\widetilde\rho=1-\varrho g(\varphi)+O(\varrho^2),
\qquad
\widetilde\varphi=f(\varphi)+O(\varrho),
\qquad
\theta=\varrho h(\varphi)+O(\varrho^2),
\]
gauge transforming the connection in \eqref{Ads} and restricting to the boundary, gives
\[
\iota^*B=
\begin{pmatrix}
\dfrac{g'+hf'}{2g} & \dfrac{f'}{g}\\[6pt]
* &-\dfrac{g'+hf'}{2g}
\end{pmatrix}d\varphi.
\]
Thus the DS conditions are
\[
\frac{f'}g=1,
\qquad
g'+hf'=0,
\]
or equivalently
\[
g=f',
\qquad
h=-\frac{f''}{f'}.
\]
Since \(f'>0\), these conditions can be achieved by a boundary-fixing diffeomorphism adjusting the normal derivative and a coframe rotation adjusting its first normal derivative. Both transformations can be supported in a boundary collar, leaving a neighbourhood of the marked point unchanged, hence the DS slice exists.

With this normalization, the remaining coefficient is
\begin{equation}\label{formula_for_T}
\iota^*B=
\begin{pmatrix}
0&1\\
-T_f&0
\end{pmatrix}d\varphi,
\qquad
T_f=\frac{\beta^2}{4}(f')^2+\frac12\mathcal{S}(f),
\end{equation}
with $\mathcal{S}(f)$ the Schwarzian derivative of $f$, as in \eqref{schwarzian}.
 
\begin{proposition}\label{dsform}
    On the DS slice, $\omega_{\rm out}$ acquires the form:
    $$
\omega_{\rm out}|_{\rm DS} = \int_{S^1}\left( \frac{1}{4} \delta \log(f')) \wedge \delta(\log(f')') - \frac{\beta^2}{4} \delta f \wedge \delta f'\right) \, d\phi.
    $$
\end{proposition}

\begin{proof}
 On the DS slice, \(g=f'\), so \(f'/g=1\) and consequently
\[
\delta\!\left(\frac{f'}g\right)=0.
\]
The mixed term in \eqref{w1} therefore vanishes, and substitution gives
\[
\omega_{\rm out}\big|_{\mathrm{DS}}
=
\int_{S^1}
\left(
\frac14\,\delta\log f'\wedge
\partial_\varphi(\delta\log f')
-\frac{\beta^2}{4}\,\delta f\wedge\delta f'
\right)d\varphi.
\]
\end{proof}
Combining the preceding results, we obtain a canonical closed \(2\)-form on Teichmüller space and identify its relation to the Virasoro coadjoint orbits. The following theorem gives an analogue of Theorem~\ref{AMtheorem} for general cone parameter $\beta$, extending the discussion from the first exceptional orbit to elliptic and exceptional orbits.

\begin{theorem}
    Fix \(\beta>0\).

(i) For each class \([g]\in\operatorname{Teich}^{\beta}(\mathbb D,p)\), one can choose a representative \(g\) and an adapted coframe such that the associated flat connection is in Drinfeld–Sokolov form near the ideal boundary.

(ii) The 2-form $\omega_{\mathrm{AB}}-\omega_{\rm in}$, restricted to the Drinfeld–Sokolov slice, descends to a canonical closed two-form \(\omega_\beta\) on \(\operatorname{Teich}^{\beta}(\mathbb D,p)\).

(iii) The natural action of \(\operatorname{Diff}^{+}(S^1)\) is a Virasoro action with moment map
\[
\mu_\beta([g])=\partial_\varphi^2+T(\varphi),
\]
where \(T(\varphi)\) is the coefficient of the Drinfeld–Sokolov connection. In the boundary parametrization by \(f\), this coefficient is
\[
T_f=\frac{\beta^2}{4}(f')^2+\frac12 S(f).
\]

(iv) If \(\beta\notin\mathbb Z_{>0}\), the moment map \(\mu_\beta\) is a bijection onto $\mathcal{O}_{\beta^2/4}$, the form \(\omega_\beta\) is symplectic and it coincides with the KKS symplectic form on the elliptic Virasoro coadjoint orbit corresponding to $T=\frac{\beta^2}{4}$.

If \(\beta=n\in\mathbb Z_{>0}\) the moment map onto its image identifies with the natural projection
\[
\pi_n:
\operatorname{Teich}^{n}(\mathbb D,p)
\cong S^1\backslash\operatorname{Diff}^{+}(S^1)
\longrightarrow
\operatorname{PSU}^{(n)}(1,1)
\backslash\operatorname{Diff}^{+}(S^1)
=\mathcal O^{(n)}.
\]
Moreover, the form $\omega_n$ is presymplectic and it coincides with the pull back of the KKS form corresponding to $T(\varphi)=\frac{n^2}{4}$ under $\pi_n$ 
\[
\omega_n=\pi_n^*\omega_{\mathrm{KKS}}^{(n)}.
\]
\end{theorem}

\begin{proof}
The proof follows by putting together the previous results. Recall that
    every Teichmüller class admits a DS representative, since the normalization can be performed in a boundary collar, leaving the cone point unchanged. By Theorem \ref{2form-tip-boundary} and Proposition \ref{dsform},
\[
(\omega_{\mathrm{AB}}-\omega_{\rm in})|_{\mathrm{DS}}
=\omega_{\rm out}|_{\mathrm{DS}}
\]
is closed and depends only on the boundary diffeomorphism \(f\). It is therefore independent of the interior extension and the coframe, and is basic under boundary-fixing diffeomorphisms. The remaining ambiguity \(f\mapsto f+c\) preserves this expression, and its infinitesimal generator \(V\), characterized by \(\delta f(V)=1\), satisfies
\[
\iota_V(\omega_{\rm out}|_{\mathrm{DS}})
=-\frac{\beta^2}{4}\int_{S^1}\delta f'\,d\varphi=0.
\]
Hence the form descends canonically to
\[
\operatorname{Teich}^{\beta}(\mathbb D)
\cong S^1\backslash\operatorname{Diff}^{+}(S^1).
\]
Denote the descended form by \(\omega_\beta\), and let
\[
q:\operatorname{Diff}^+(S^1)
\longrightarrow S^1\backslash\operatorname{Diff}^+(S^1)
\]
be the quotient map. Comparing Proposition \ref{dsform} with \eqref{eq: KKS},
for \(T=\beta^2/4\), gives
\[
q^*\omega_\beta=\nu_T^*\omega_T.
\]
By Equation \eqref{formula_for_T}, the boundary potential $T_f$ is precisely the coadjoint transform of the constant potential $T_\beta=\beta^2/4$ under the action $\eqref{eq:Virasoro action}$. Hence $\mu_\beta \circ q=\nu_{T_\beta}$.

If \(\beta\notin\mathbb Z_{>0}\), the stabilizer of \(T\) is \(S^1\).
Hence the identification of Teichmüller space with the corresponding
elliptic orbit identifies \(\omega_\beta\) with its KKS symplectic form.

If \(\beta=n\in\mathbb Z_{>0}\), the stabilizer is
\(\operatorname{PSU}^{(n)}(1,1)\), and
\(\nu_T=\pi_n\circ q\). Thus
\[
q^*\omega_n=q^*\pi_n^*\omega_{\mathrm{KKS}}^{(n)}.
\]
It follows that
\[
\omega_n=\pi_n^*\omega_{\mathrm{KKS}}^{(n)},
\]
and hence $\omega_n$ 
is presymplectic.
\end{proof}

\appendix
\section{Hyperbolic cone metrics on the disk}
\label{app-conemetrics}

In this appendix, we recall the relevant results on simple singularities
from the work of  Hulin and Troyanov \cite{hulin-troyanov, troyanov}, as well as McOwen  \cite{mcowen}, and give a detailed proof
of Theorem \ref{thm:Troyanov}. Existence follows from uniformization argument. Uniqueness is established by the Schwarz lemma,
using the behavior of the conformal ratio at the marked point and
at the ideal boundary.

Throughout this appendix, singularity orders are understood in the sense of Definition \ref{simple_singularity}. Our order $\beta>0$ corresponds to order $\alpha=\beta-1$ in the convention of \cite {troyanov, hulin-troyanov}.

\subsection{Simple singularities and uniqueness}
A metric with a simple singularity of order $\beta$, has cone angle $2\pi\beta$ at $p$, it represents the divisor
$(\beta-1)p$ and has the local form
\[
g=e^{2u}|z|^{2(\beta-1)}|dz|^2,
\qquad
u\in L^1_{\mathrm{loc}},
\qquad
\Delta_0u\in L^1_{\mathrm{loc}}.
\]
Here $z(p)=0$, the integrability conditions hold across $p$, and
$\Delta_0=\partial_x^2+\partial_y^2$.

\begin{theorem}[McOwen \cite{mcowen}]
Let $g=e^{2v}|dz|^2$ be a conformal metric on the punctured unit disk,
smooth away from the origin. If its curvature satisfies
\[
-b|z|^{\ell}\leq K(z)\leq-a|z|^{\ell}
\]
for some $\ell\in\mathbb R$ and $a,b>0$, then the origin is a simple
singularity of $g$.
\end{theorem}

For $K=-1$, the hypotheses hold with $\ell=0$ and $a=b=1$.
This criterion guarantees that the singularity is simple but does
not determine its order. In Theorem \ref{thm:Troyanov}, we fix $\beta>0$ in advance
and require the metric to have a simple singularity of order
$\beta$ at the marked point.

\begin{proposition}[Schwarz lemma for the disk; Hulin--Troyanov
{\normalfont[\cite{hulin-troyanov}, Theorem~3.1]}]\label{schwarz_lemma}
Let $g_1$ and $g_2$ be conformal metrics of curvature $-1$ on
$\mathbb{D}\setminus\{p\}$, each with a simple singularity of order
$\beta$ at $p$, where $\beta>0$. Write $g_2=e^{2v}g_1$.
If $g_1$ is complete at the infinity boundary of $\mathbb{D}$ and $v$ extends
continuously across $p$, then $g_2\leq g_1$ on $\mathbb{D}\setminus\{p\}$.
\end{proposition}

The proof of the Schwarz lemma relies on the generalized maximum
principle of Yau; see \cite{troyanov} and \cite{hulin-troyanov}. The next lemma
provides the control near the marked point needed to apply
Proposition \ref{schwarz_lemma} in the uniqueness proof.

\begin{lemma}
\label{lem:cone-local-regularity}
Suppose
\[
g=e^{2u}|z|^{2(\beta-1)}|dz|^2
\]
has curvature $-1$ and a simple singularity with $\beta>0$ at $z=0$.
Then $u$ has a continuous representative across the origin.
Consequently, two such metrics with the same parameter $\beta$ have a
bounded logarithmic  ratio near the origin.
\end{lemma}

\begin{proof}
Away from the origin, the curvature equation gives
\[
\Delta_0u=|z|^{2(\beta-1)}e^{2u}\geq0.
\]
Since $\Delta_0u$ is an $L^1$ function on the whole coordinate disk,
this identity holds distributionally across the origin as well.
Hence $u$ is subharmonic and locally bounded above.
The right-hand side belongs locally to $L^q$ for some $q>1$:
if $0<\beta<1$, choose
\[
1<q<(1-\beta)^{-1},
\]
while for $\beta\geq1$ it is locally bounded.
Interior elliptic regularity gives
\[
u\in W^{2,q}_{\mathrm{loc}}(U),
\]
where $U$ is an open neighbourhood in the disk containing the origin.
By Sobolev embedding in dimension two, $u$ therefore admits
a continuous representative on $U$, including the origin.

For two such metrics with the same parameter $\beta$, write
\[
g_i=e^{2u_i}|z|^{2(\beta-1)}|dz|^2.
\]
The common singular power cancels, giving
\[
g_2=e^{2v}g_1,
\qquad v=u_2-u_1.
\]
Since both $u_1$ and $u_2$ extend continuously across the origin,
so does $v$, with
\[
\lim_{z\to0}v(z)=u_2(0)-u_1(0)\in\mathbb R.
\]
Consequently, 
$v$ is bounded
in a neighborhood of the origin.
\end{proof}
\subsection{Proof of Theorem \ref{thm:Troyanov}}
\begin{proof}
Fix a smooth complex structure $I$ on $\overline{\mathbb D}$,
a point $p\in\mathbb D$, and a parameter $\beta>0$.
We prove that there is a unique compatible hyperbolic $0$-metric with a simple singularitiy of order $\beta$ at $p$.

\emph{Existence.}
Uniformization, together with smooth boundary regularity, gives a
biholomorphism of the interiors
\[
w:(\mathbb D,I)\longrightarrow(\mathbb D,J),
\qquad
w(p)=0,
\]
extending to a smooth diffeomorphism of the closed disks.
Consider the model
\[
g_{\mathbb D}^{\beta}
=
\frac{4\beta^2|z|^{2(\beta-1)}}
     {(1-|z|^{2\beta})^2}|dz|^2,
\qquad
g=w^*g_{\mathbb D}^{\beta}.
\]
The model has curvature $-1$ and cone angle $2\pi\beta$ at the origin.
Since $w$ is a local diffeomorphism at $p$, the pullback has the same
cone angle there.
On a collar of the outer boundary, put $\varrho=1-|z|$.
Then
\[
\varrho^2g_{\mathbb D}^{\beta}
=
\left(
\frac{2\beta|z|^{\beta-1}(1-|z|)}
     {1-|z|^{2\beta}}
\right)^2|dz|^2
\]
extends smoothly and positively to $|z|=1$.
Thus $g_{\mathbb D}^{\beta}$, and hence its pullback, satisfies the
$0$-metric condition.
The normalized map $w$ is unique up to a rotation of the target, and
rotations preserve $g_{\mathbb D}^{\beta}$, so the construction is
independent of this choice.

\emph{Uniqueness.}
Let $g_1$ and $g_2$ 
be two hyperbolic 0-metrics having a simple singularity of order $\beta$ at $p$. 
In a common holomorphic coordinate centred at $p$, write
\[
    g_i=e^{2u_i}|z|^{2(\beta-1)}|dz|^2, \qquad i=1,2.
\]
By Lemma \ref{lem:cone-local-regularity}, each $u_i$ extends continuously across $p$. Consequently,
the logarithmic conformal ratio defined by
\[
    g_2=e^{2v}g_1,\qquad v=u_2-u_1\quad\text{near }p,
\]
extends to a continuous function on $\mathbb{D}$, smooth away from $p$.

Both metrics are complete at infinity. Indeed, near the ideal boundary
they have the form
\[
    g_i=\varrho^{-2}\bar g_i,
\]
where $\varrho$ is a boundary defining function and $\bar g_i$ extends
smoothly and positively to the boundary. 
Since $K_{g_1}=K_{g_2}=-1$, the hypotheses of the Schwarz lemma for
metrics with singularities of the same order are satisfied.
Applying it with $(g,h)=(g_1,g_2)$ gives $g_2\leq g_1$.
Interchanging the metrics gives $g_1\leq g_2$. Thus $g_1=g_2$.
\end{proof}

For each fixed $\beta>0$, Theorem \ref{thm:Troyanov} thus identifies a marked conformal
disk with its unique compatible hyperbolic cone metric.
Uniqueness makes this correspondence natural under pullback by
diffeomorphisms.
This is used in Proposition \ref{hypbeta} and in the subsequent
description of Teichm\"uller space.

\end{document}